\documentclass[12pt,reqno,NoRunningHeads]{amsart}
\usepackage{latexsym}
\usepackage{amsmath}
\usepackage{amssymb}
\usepackage{amsfonts}
\usepackage{verbatim}
\usepackage[latin1]{inputenc}
\usepackage{color}

\newtheorem{Th}{Theorem}[section]
\newtheorem{assumption}{Assumption}[section]
\newtheorem{corollary}[Th]{Corollary}

\newtheorem{proposition}[Th]{Proposition}

\newtheorem{Ex}[Th]{Example}
\newtheorem{Def}[Th]{Definition}

\newtheorem{rk}[Th]{Remark}

\let \ssection=\section
\renewcommand{\section}{\setcounter{equation}{0}\ssection}

\def\^#1{\if#1i{\accent"5E\i}\else{\accent"5E #1}\fi}
\def\"#1{\if#1i{\accent"7F\i}\else{\accent"7F #1}\fi}

\newcommand{\RR}{\mathbb R}

\newcommand{\ZZ}{\mathbb Z}

\newcommand{\<}{\langle}
\renewcommand{\>}{\rangle}
\def\^#1{\if#1i{\accent"5E\i}\else{\accent"5E #1}\fi}
\def\"#1{\if#1i{\accent"7F\i}\else{\accent"7F #1}\fi}
\allowdisplaybreaks

\begin{document}
\title[Accuracy of  Space-Time Approximations]
{Rate of Convergence of Space Time Approximations\\
for stochastic evolution equations}
\author[I. Gy\"ongy]
{Istv\'an Gy\"ongy}
\thanks{This paper was written 
while the first named author was visiting
the University of Paris 1. 
The research of this author is partially supported by
EU Network HARP}

 \author[A. Millet]{Annie Millet}
\thanks{The research of the second named author is partially supported
by the research project BMF2003-01345}



\begin{abstract}
Stochastic evolution equations
in Banach spaces
with unbounded nonlinear drift and diffusion operators 
driven by a finite dimensional
Brownian motion are
considered. Under some regularity 
condition assumed for the solution,
the rates of
convergence of various numerical approximations are
estimated under strong monotonicity and Lipschitz
conditions. The abstract setting involves general
consistency conditions and is then applied to a class of
quasilinear stochastic PDEs of parabolic type.
\end{abstract}

\maketitle

\section{Introduction}
                                            \label{intro}
Let
$V\hookrightarrow  H\hookrightarrow   V^*$
be a
{\it normal triple} of spaces with dense and
continuous embeddings, where $V$
is a separable and reflexive Banach space,
$H$ is a Hilbert space,
identified with its dual 
 by
means 
of the inner product in $H$,
and $V^*$ is the dual of $V$.

 Let
$W=\{W(t):t\geq 0\}$ be a
$d_1$-dimensional
Brownian motion carried by
a stochastic basis
$(\Omega, \mathcal F, ({\mathcal F}_t)_{t\geq 0}, P)$. 
Consider the stochastic evolution equation
\begin{equation}                                      \label{u}
 u(t)=u_0 + \int_0^t A(s,u(s))\, ds
+  \sum_{k=1}^{d_1} \int_0^t
 B_k(s,u(s))\, dW^k(s)\,, \quad t\in[0,T] 
\end{equation}
in the triple $V\hookrightarrow  H\hookrightarrow   V^*$, 
with a given $H$-valued $\mathcal F_0$-measurable 
random variable $u_0$, and given operators 
$A$ and $B=(B_k)$, mapping $[0,\infty)\times\Omega\times V$ 
into $V^{\ast}$ and $H^{d_1}:=H\times\dots\times H$, 
respectively. Let $\mathcal P$ denote the 
$\sigma$-algebra of the predictable subsets of 
$[0,\infty)\times\Omega$, and let $\mathcal B(V)$, 
$\mathcal B(H)$ and $\mathcal B(V^{\ast})$ be 
the Borel $\sigma$-algebras of $V$, $H$ and $V^{\ast}$, 
respectively. Assume that $A$ and $B_k$ are 
${\mathcal P}\otimes {\mathcal B}(V)$-measurable 
with respect to the $\sigma$-algebras 
$\mathcal B(V^{\ast})$ and $\mathcal B(H)$, respectively. 

It is well-known that for any $T>0$ 
equation \eqref{u} admits a unique 
solution $u$ if $A$ is {\it hemicontinuous} in $v\in V$, and 
$(A,B)$ satisfies a {\it monotonicity}, {\it coercivity} 
and a {\it linear growth} condition (see \cite{KR}, 
\cite{Pa} and \cite{R}). 
In \cite{GyMi} it is shown that under these conditions the
solutions of various implicit and explicit schemes converge
to $u$. In \cite{GyMi2} the rate of convergence of implicit
Euler approximations is estimated under more
restrictive hypotheses:
$A$ and $B$ satisfy a strong monotonicity
condition, 
$A$ is Lipschitz continuous in $v\in V$, 
and the solution $u$
satisfies some  regularity conditions. 
Then Theorem 3.4
 from \cite{GyMi2} in the case of 
time independent operators
$A$ and
$B$ reads as follows.
For the implicit Euler approximation
$u^{\tau}$, corresponding  to the mesh size $\tau=T/m$ of
the partition of $[0,T]$, one has
\[
E\max_{i\leq  m}|u(i\tau)-u^{\tau}(i\tau)|_H^2
+ \tau E\sum_{i\leq  m}\|u(i\tau)-u^{\tau}(i\tau)\|_V^2
\leq C\tau^{2\nu}, 
\]
where $C$ is a constant, 
independent of $\tau$, and $\nu\in
]0,\frac{1}{2}]$ is a constant from the regularity 
condition imposed on $u$. 

In this paper, we study space and space-time
approximations schemes 
 for equation \eqref{u} in a general framework. 
In order to obtain rate of convergence estimates 
we need to require more regularity from the solution $u$ 
of equation \eqref{u} than what we can express in terms 
of the spaces $V$ and $H$. Therefore in our setup 
we introduce additional  
Hilbert spaces $\mathcal V$ and $\mathcal H$ 
such that 
$\mathcal V\hookrightarrow \mathcal H\hookrightarrow V$, 
where $\hookrightarrow$ denotes continuous embeddings. 
In examples these are Sobolev spaces such that 
$\mathcal H$ and $\mathcal V$ satisfy stronger 
differentiability conditions than $V$ and $\mathcal H$, 
respectively. Our regularity conditions 
on the solution $u$ are introduced in section 2 and 
labeled as {\bf(R1)} and 
{\bf(R2)}. 
In connection with these, we introduce also 
condition {\bf(R3)}, requiring more regularity 
from $A$ and $B$. 
Furthermore,  condition {\bf(R4)} on H\"older 
continuity in time of $A$ and $B$ is needed  
for schemes involving time discretization. 
We collect these conditions 
in Assumption \ref{assumption R} and call them 
{\it regularity conditions}. 

In order to formulate `space discretizations',  
we consider for any integer $n\geq1$ a normal triple
\begin{equation}                              \label{triplen}
V_n \hookrightarrow  H_n  \hookrightarrow   V_n^*, 
\end{equation} 
 the `discrete' counterpart of 
$V \hookrightarrow  H  \hookrightarrow   V^* $, 
and a bounded linear operator
 $$
\Pi_n : V\to V_n, 
$$
connecting $V$ 
to $V_n$.  
We have in mind 
discrete Sobolev spaces, wavelets and 
finite elements spaces,  as examples 
 for $V_n$.

The space discretization scheme  
for 
equation \eqref{u} is a stochastic 
evolutional equation 
in the triple \eqref{triplen}. 
We define it by replacing the 
operators $A$, $B$ and the initial value $u_0$ 
in 
equation \eqref{u} by some  
$\mathcal P\otimes \mathcal B(V_n)$-measurable 
operators 
$$
A^n: [0,\infty[\times \Omega\times V_n
\to V_n^*,
\quad
B^n: [0,\infty[\times \Omega\times V_n\to
H_n^{d_1} 
$$
and by 
an $H_n$-valued $\mathcal F_0$-measurable random variable 
$u^n_0$, 
respectively, such that $A^n$ and $B^n$ satisfy in the 
triple \eqref{triplen}
the strong monotonicity condition, the linear growth condition, 
$A^n$ is hemicontinuous and $B^n$ is Lipschitz continuous 
in $v\in V_n$.  
 These are the conditions 
{\bf (S1)}--{\bf (S4)} in Assumption 3.1, which imply, in
particular, the  existence and uniqueness of a solution $u^n$ to
our scheme.     We relate   
$A_n$ and
$B_n$ to $A$ and $B$ 
via a {\it consistency condition}, {\bf (Cn)} below.
Then
assuming {\bf (S1)}--{\bf (S4)}, under the regularity and consistency
conditions {\bf (R1)}, {\bf (R3)} 
and {\bf (Cn)}
we have
$$
E \sup_{0\leq t\leq T} 
|\Pi_n u(t) - u^n(t)|_{H_n}^2
+E\int_0^T \| \Pi_n u(t)- u^n(t)\|_{V_n}^2
dt 
\leq CE|\Pi_n u_0 - u^n_0|_{H_n}^2+
C\varepsilon_n^2,
$$
where 
 $C$ is a constant, independent 
of $n$, and $\varepsilon_n>0$
 is a constant from 
{\bf(Cn)}. This is Theorem \ref{raten} below, our main 
result on the accuracy of approximations 
by space discretizations.

For an integer $m\geq 1$  
we consider the grid 
$\{t_i=i\, \tau : 0\leq i\leq m\}$ with mesh-size 
$\tau=T/m$. 
We define on this grid the space-time
implicit  and the
 space-time explicit
approximations, $\{u^{n,\tau}_i\}_{i=0}^m$ 
and $\{u^n_{\tau,i}\}_{i=0}^m$,   
 respectively, 
by
\begin{eqnarray*}
u^{n,\tau}_{i+1}&=& u^{n,\tau}_{i} + \tau
A^{n,\tau}_{i+1}\big(u^{n,\tau}_{i+1}\big) +
\sum_k B^{n,\tau}_{k,i}\big( u^{n,\tau}_{i}\big)
\, \big( W^k(t_{i+1})- W^k(t_{i})\big),\\
u^n_{\tau,i+1} & =& u^n_{\tau,i} + \tau\,
A^{n,\tau}_i(u^n_{\tau,i}) + \sum_k
B^{n,\tau}_{k,i}(u^n_{\tau,i}) \big( W^k(t_{i+1})-W_k(t_i) \big),
\end{eqnarray*}
for $i=0,\dots,m-1$ with 
some 
$V_n$-valued $\mathcal F_0$-measurable random variables 
$u_0^{n,\tau}$ and $u_{0,\tau}^{n}$, and 
with 
some $\mathcal F_{t_i}\otimes\mathcal B(V)$-measurable  
operators
$$
A^{n,\tau}_{i}:\Omega\times V_n\to V_n^*,
\quad
B^{n,\tau}_{k,i}: \Omega\times V_n\to
H_n^{d_1},   
$$
 such that $A^{n,\tau}_{i}$,  
$B^{n,\tau}_{k,i}$ 
satisfy strong monotonicity and linear 
growth conditions  and $A^{n,\tau}_{i}$ is Lipschitz 
continuous in $v\in V_n$. These conditions,  listed as 
{\bf(ST1)}--{\bf(ST3)} in Assumption 4.1  below, 
are similar 
to conditions {\bf(S1)}--{\bf(S3)}, 
except that instead of the hemicontinuity, 
the much stronger assumption of Lipschitz continuity 
is assumed on $A^{n,\tau}_i$.  
The operators $A^{n,\tau}_{i}$ and   
$B^{n,\tau}_{k,i}$ 
are related 
to $A$ and $B$ by a consistency condition 
{\bf(Cn$\tau$)}  stated  below.
Then if 
$\sup_{n,m}E|u_0^{n,\tau}|_{H_n}^2<\infty$ 
and equation \eqref{u} satisfies the regularity 
conditions {\bf(R1)}--{\bf(R4)} 
from Assumption \ref{assumption R}, we
have the estimate
$$
 E \sup_{0\leq i\leq m}
|\Pi_n u(t_i) - u^{n,\tau}_i|_{H_n}^2 
+E\sum_{0\leq i\leq m}
\|\Pi_n u(t_i)- u^{n,\tau}_i\|_{V_n}^2 \tau
$$
$$
 \leq 
 CE|\Pi_n u_0 - u^{n,\tau}_0|_{H_n}^2 
+C( \tau^{2\nu} +\varepsilon_n^2), 
$$
with a constant $C$, 
independent of $n$ and $\tau$,  where 
$\nu\in
]0,\frac{1}{2}]$ is the H\"older exponent 
from condition {\bf(R4)} 
on the regularity of the operators $A$ and $B$ in 
time, and $\varepsilon_n$ is 
from {\bf(Cn$\tau$)}. 
This is Theorem \ref{th1spi}, our main result 
on implicit space-time approximations. 
In our main result, Theorem \ref{thexp}, on the explicit
space-time approximations  we have the same estimate for
$u^{n}_{\tau,i}$  in place of $u^{n,\tau}_{i}$ if,  in addition 
to the conditions of Theorem \ref{th1spi}, 
as in \cite{GyMi},
a {\it stability relation} 
between the time mesh $\tau$ and a
space approximation parameter is satisfied.

Finally,  we present as examples  a class of quasi-linear
stochastic partial differential equations (SPDEs) 
and linear SPDEs of parabolic
type. We show that they satisfy the conditions 
of the abstract results, Theorems \ref{raten}, 
\ref{th1spi} and \ref{thexp}, when we use wavelets, or finite 
differences.  In particular, we
obtain rate of convergence results for space and space-time
approximations of linear parabolic SPDEs, 
among them 
for the Zakai equation of nonlinear filtering. 
We would like to mention that as far as we know, discrete 
Sobolev spaces are applied first in \cite{HY} to space
discretizations and explicit  space-time discretizations 
of linear SPDEs, and it inspired our approach to finite  
difference schemes. Our abstract results can also be applied 
to finite elements approximations. To keep down the size of 
the paper we will consider such applications elsewhere.

We denote by $K$, $L$,  $M$ and $r$  some fixed 
constant, and by $C$  
 some constants
which, as usual,  
can change from line to line. 
For given constants $a\in \RR^k$ the notation 
$C(a)$ means that the constant depends on
$a$.
 Finally, when $(X, |\cdot |_X)$ 
and $(Y,|\cdot|_Y)$ denote two Banach
spaces such that $X$ is continuously 
embedded in $Y$, given $y\in Y$ the inequality
$|y|_X<+\infty$ means that $y\in X$.

\section{Conditions on equation \eqref{u} and 
on the  approximation spaces}
\subsection{Conditions on equation \eqref{u}}

Let $(\Omega, {\mathcal  F}, ({\mathcal F}_t)_{t\geq 0}, P)$ be a
stochastic basis, satisfying the usual conditions, i.e.,
$({\mathcal F}_t)_{t\geq0}$ is an increasing right-continuous
family of sub-$\sigma$-algebras of $\mathcal F$ such that
${\mathcal F}_0$ contains every $P$-null set. Let
$W=\{W(t): t\geq0\}$ be a $d_1$-dimensional
Wiener martingale with respect to
$({\mathcal F_t})_{t\geq 0}$, i.e., $W$ is an
${\mathcal F}_t$-adapted Wiener process
with values in ${\mathbb R}^{d_1}$ such that
$W(t)-W(s)$ is independent of ${\mathcal F}_s$ for all
$0\leq s\leq t$. 
We use the notation 
${\mathcal P}$ for the sigma-algebra 
of predictable subsets of $[0,\infty)\times\Omega$. 
If $V$ is a Banach space then ${\mathcal B}(V)$ 
denotes the sigma-algebra generated by the (closed) 
balls in $V$.

Let $V$ be a separable 
reflexive Banach space 
embedded densely and 
continuously into a Hilbert space $H$,  
which is 
 identified with its dual $H^{\ast}$ 
by means of the inner product 
$(\cdot,\cdot)$ in $H$. Thus we have 
a {\it normal triple}
$$
V\hookrightarrow H\hookrightarrow V^*, 
$$
where $H\hookrightarrow V^*$ is the adjoint 
of the embedding $V\hookrightarrow H$. 
Thus $\<v,h\>=(v,h)$ for all
$v\in V$ and $h\in H^*=H$,
where $\<v,v^*\>=\<v^*,v\>$
denotes the duality product
of $v\in V$, $v^*\in V^*$,
and $(h_1,h_2)$ denotes the
inner product of $h_1,h_2\in H$. 
We assume, without loss of generality, that 
$|v|_H\leq \|v\|_H$ for all $v\in V$, where 
$|\cdot|_H$ and $\|\cdot\|_V$ denote the norms 
in $H$ and $V$, respectively. For elements 
$u$ from a normed space $\mathbb U$ the notation 
$|u|_{\mathbb U}$ means the norm of 
$u$ in $\mathbb U$.   
\medskip
  
Let $A$ and $B=(B_k)_{k=1}^{d_1}$ 
be $\mathcal
P\otimes\mathcal B(V)$-measurable  mappings 
from 
$
[0,\infty)\times\Omega\times V
$
into  $V^{\ast}$ and $H^{d_1}$, 
respectively.   
Given an $H$-valued ${\mathcal F}_0$-measurable random variable
$u_0$ 
consider the initial value problem
\begin{equation} \label{eq2.1}
du(t)=A(t,u(t))\,dt+\sum_kB_k(t,u(t))\,dW^k(t), 
\quad u(0)=u_0
\end{equation}                                        \label{v}
on a fixed time interval $[0,T]$.

\begin{assumption}                              \label{assumption 0}
The operators $A$ and $B$ satisfy the following conditions. 

{\bf (i)} (Monotonicity of $(A,B)$)\quad 
Almost surely for all $t\in [0,T]$ and  $u,v\in V$, 
\begin{equation*}                                  \label{1}
  2\<u-v,A(t,u)-A(t,v)\>
+ \sum_k |B_k(t,u)-B_k(t,v)|^2_H 
\leq K|u-v|^2_H,   
\end{equation*}

{\bf (ii)} (Coercivity of $(A,B)$)
Almost surely for all $t\in [0,T]$ and  $u,v\in V$,
\begin{equation}                                  \label{2}
  2\<u,A(t,u)\>
+ \sum_k |B_k(t,u)|^2_H 
 +\mu \|u\|^2_V 
\leq K|u|^2_H + f(t),   
\end{equation}

{\bf (iii)} (Growth conditions on $A$ and $B$) \quad 
Almost surely for all $t\in [0,T]$ and 
$u\in V$, 
\begin{equation*}                                 \label{3}
|A(t,u)|_{V^*}^2 
\leq  K_1 \|u\|^2_V  +f(t),
\quad \sum_k |B_k(t,u)|_H^2 \leq K_2
\|u\|^2_V + f(t),    
\end{equation*}

{\bf(iv)} (Hemicontinuity of $A$) \quad 
Almost surely for all $t\in [0,T]$ and  $u,v,w\in V$,
\begin{equation}                              \label{4}
   \lim_{\varepsilon\to0}\<w,A(t,u+\varepsilon v)\>=
\<w,A(t,u)\>, 
\end{equation} 
where $\mu>0$, $K\geq0$, $K_1\geq0$ and 
$K_2\geq0$  are some constants, and $f$ is a  nonnegative 
$(\mathcal F_t)$-adapted  stochastic process such that 
\begin{equation}                                       \label{fg}
E\int_0^Tf(t)\,dt<\infty.
\end{equation} 
\end{assumption}
The following definition of solution is classical. 
\begin{Def}                                             \label{sol}
An $H$-valued adapted continuous process 
$u=\{u(t):t \in[0,T ] \}$ is a solution to
equation \eqref{u} on $[0,T]$ if almost surely  
$u(t)\in V$ for almost every $t\in[0,T]$, 
\begin{equation*} 
\int_0^T  \|u(t)\|^2_V\, dt <\infty\, , \label{definition*}
\end{equation*} 
and 
\begin{equation*}                                  \label{solu}
( u(t), v)  =(u_0,v)
+ \int_0^t \langle A(s,u(s)),v
\rangle \, ds
+ \sum_{k}   \int_0^t (B_k(s,u(s)),v)\,
dW^k(s)
\end{equation*}
holds for all $t\in [0,T]$ and $v \in V$. 
We say that the solution to \eqref{eq2.1}  on $[0,T]$ is unique 
if for any solutions $u$ and $v$ to 
\eqref{eq2.1}  
on $[0,T]$  we have 
$$
P\Big(\sup_{t\in[0,T]}|u(t)-v(t)|_H>0\Big)=0.
$$
\end{Def}

The following result is well known,   see  
\cite{KR}, \cite{Pa}, \cite{R}.

\begin{Th}                               \label{existu}  
Let Assumption \ref{assumption 0} hold. Then 
\eqref{v} has a
unique solution $u$. 
Moreover, if $E|u_0|_H^2<\infty$, then    
\begin{align}                                      \label{bounduH}
E\sup_{t\in[0,T]} & |u(t)|^2_H 
+E\int_0^T \|  u  (s)\|^2_V\,ds \nonumber \\
& \leq CE|u_0|_H^2 
+ CE\int_0^T \big(f(t)+g(t)\big)\, dt<\infty ,
\end{align}
where $C$ is a constant  
depending only on the constants $\lambda$, 
$K$ and $K_2$. 
\end{Th}

If Assumption \ref{assumption 0} is satisfied 
then one can also show the convergence 
of approximations, obtained by 
 various discretization schemes,  
to the solution $u$ (see \cite{GyMi}).  
To estimate the rate of convergence of implicit 
time discretization schemes the following stronger 
assumption on $A$ and $B$ are used in \cite{GyMi2}

\begin{assumption}                          \label{E}
The operators $A$, $B$ satisfy the following conditions 
almost surely. 

{\bf (1)} (Strong monotonicity) 
For all $t\in [0,T]$, $u,v\in V$, 
$$
  2\<u-v,A(t,u)-A(t,v)\>
+ \sum_k |B_k(t,u)-B_k(t,v)|_H^2 
$$
\begin{equation*}                                         \label{R01}
\leq -\lambda \|u-v\|^2_V 
+L|u-v|_H^2,   
\end{equation*}

{\bf (2)} (Growth conditions on $A$ and $B$) 
For all $t\in [0,T]$,
$u\in V$, 
\begin{equation}                                          \label{R02}
|A(t,u)|_{V^*}^2 
\leq  K_1 \|u\|^2_V  +f(t),\quad 
\sum_k |B_k(t,u)|_H^2 \leq K_2 \|u\|^2_V
+ g(t).   
\end{equation}

{\bf(3)} (Lipschitz condition on $A$) \quad 
For all $t\in [0,T]$, $u,v\in V$, 
\begin{equation}                                            \label{R03}
   |A(t,u)- A(t,v)|_{V^*}^2\leq L_1   \|u-v\|^2_V,  
\end{equation}
where $\lambda>0$, $K\geq0$, $K_1\geq0$, $K_2\geq0$ 
are constants, and $f$ and $g$ are non-negative 
adapted processes satisfying \eqref{fg}
\end{assumption}

\begin{rk}                               \label{remark 28.09.08}
It is easy to see that due to 
{\bf(1)}--{\bf(2)}, the coercivity condition 
\eqref{2} holds 
with 
$\mu=\lambda/2$ and a constant $K=K(\lambda,L,K_2)$.   
\end{rk}
\begin{rk}                               \label{remark 29.09.08}
It is easy to show that 
{\bf(1)} and  {\bf(3)} imply that $B=(B_k)$ 
is Lipschitz continuous in $u\in V$, i.e., 
almost surely 
\begin{equation}                           \label{29.09.08}
\sum_k|B_k(t,u)-B_k(t,v)|_H^2\leq L_2\|u-v\|^2_V
\quad \text{for all $u,v\in V$, $t\in[0,T]$}
\end{equation}
where $L_2$ is a constant depending on 
$\lambda$, $L$ and $L_1$.    
\end{rk}
In order to prove rate of convergence 
estimates for the approximation 
schemes presented in this paper,  
we need to impose additional  {\it regularity 
conditions} on equation \eqref{v} and 
on the solution $u$. Therefore we assume 
that there exist some  
separable Hilbert spaces $\mathcal V$ and 
$\mathcal H$ such that 
$$
\mathcal V\hookrightarrow \mathcal H\hookrightarrow V, 
$$ 
where $\hookrightarrow$ means continuous embedding, 
and  
introduce the following conditions. 

Let $K$,  $M$ denote some constants, 
fixed  throughout  the paper.  
\begin{assumption}                        \label{assumption R} 
(Regularity conditions)

{\bf(R1)} There is 
a unique solution $u$ of 
\eqref{v}, it takes values in 
$\mathcal V$ for $dt\times P$-almost every 
$(t,\omega)\in[0,T]\times\Omega$, $u_0\in V$ and  
\begin{equation}                                 \label{r1}
E\|u_0\|^2_V<\infty, 
\quad E\int_0^T|u(t)|^2_{\mathcal
V}\,dt=:r_1<\infty.  
\end{equation}

{\bf(R2)} There is a 
unique solution $u$ of 
\eqref{v}, it has an $\mathcal H$-valued 
stochastic modification, denoted also by $u$, 
such that 
\begin{equation*}                                  \label{r2}
\sup_{t\in[0,T]}E|u(t)|_{\mathcal H}^2=:r_2<\infty.
\end{equation*}

{\bf (R3)} Almost surely $A(t,v)\in V$,
$B_k(t,u)\in V$ and 
\begin{equation}                                    \label{R3}
\|A(t,v)\|^2_{V}\leq K|v|^2_{\mathcal V}+\xi(t), 
\quad 
\sum_k\|B_k(t,u)\|_{V}^2
\leq K|u|^2_{\mathcal H}+\eta(t) 
\end{equation}
for all $t\in[0,T]$ $v\in \mathcal V$ 
and $u\in \mathcal H $, 
where 
$\xi$ and $\eta$ are 
non-negative processes such that  
for some constant $M$ 
\begin{equation*}                                   \label{M}
 E\int_0^T\xi(t)\,dt\leq M, 
\quad 
\sup_{t\in[0,T]}E\eta(t)\leq M. 
\end{equation*}
\smallskip

{\bf (R4)} (Time regularity of $A$, $B$)   There
exists a constant
$\nu \in ]0,\frac{1}{2}]$ 
and a non-negative random variable $\eta $ 
such that $E\eta\leq M$, and almost surely 

(i)  
\begin{equation}                                       \label{HA}
\|A(s,v) - A(t,v)\|_{V}^2 
\leq 
(K\, |v|_{\mathcal V}^2 +\eta)\,|t-s|^{2\nu}
\quad
\text{for $v\in\mathcal V$},  
\end{equation}

(ii)  
\begin{equation*}                                     \label{HB}
\sum_k |B_k(s,u) -B_k(t,u)|_{V}^2 
\leq (K\, |u|_{\mathcal V}^2
+ \eta )\,  |t-s|^{2\nu} 
\quad  \text{for $u\in\mathcal V$},
\end{equation*}
for all $0\leq s<t\leq T$. 

\end{assumption}
\begin{rk}                                          \label{19.02.04}
Assume conditions {\bf(R1)} and 
{\bf(R3)} from Assumption \ref{assumption R}.     
Then the following statements hold.

(i) $u$ has a $V$-valued continuous stochastic 
modification, denoted also by $u$, such that 
$$
E\sup_{t\in[0,T]}\|u(t)\|^2\leq 3E\|u_0\|^2_V
+C(r_1+M  ); 
$$

(ii) If condition 
{\bf(R2)} from Assumption
\ref{assumption R} 
also holds, then  for $s,t\in [0,T]$, 
\begin{equation}                                  \label{Holderu}
E\|u(t)-u(s)\|^2_V \leq C|t-s|(r_1+r_2+M), 
\end{equation}
where $C$ is a constant depending
only   on $T$ and on  the constant $K$ 
from  \eqref{R3}. 
\end{rk}
\begin{proof} Define 
$$
F(t)=\int_0^tA(s,u(s))\,ds 
\quad \mbox{\rm and } \quad  
G(t)=\sum_k\int_0^tB_k(s,u(s))\,dW^k(s). 
$$
Notice that 
\begin{eqnarray*}
E\int_0^T \|A(s,u(s))\|^2\,ds\leq KE\int_0^T
|u(s)|^2_{\mathcal V}\,ds+E\int_0^T\xi(s)\,ds
&=:&M_1<\infty,  \\
\sum_k\int_0^T E\|B_k(s,u(s))\|_{V}^2\,ds
\leq K \, E\int_0^T
|u(s)|^2_{\mathcal H}\,ds 
+E\int_0^T  \eta(s)  \,ds &=:&M_2<\infty.
\end{eqnarray*} 
Hence $F$ and $G$ are   $V$-valued continuous processes,  
and by Jensen's and Doob's inequalities
$$
E\sup_{t\leq T}\|F(t)\|^2_V\leq TM_1, 
\quad
E\sup_{t\leq T}\|G(t)\|^2_{V}
\leq 4\sum_kE\int_0^T\|B_k(s,u(s))\|_{V}^2\,ds
\leq 4M_2. 
$$
Consequently, 
the process $u_0+F(t)+G(t)$ is  
a $V$-valued continuous modification of 
$u$, and statement (i) holds. Moreover,  
if  {\bf (R2)} also  holds, then 
$$
\sup_{t\in[0,T]} \sum_k E\|B_k(s,u(s))\|_{V}^2
\leq K\sup_{t\in[0,T]}
E|u(t)|^2_{\mathcal H}+\sup_{t\in[0,T]}E\eta(t):=M_3 <+\infty ,
$$
and  
\begin{align*}
E\|F(t)-F(s)\|^2_V& \leq |t-s|M_1, \\ 
E\|G(t)-G(s)\|^2_V&=
\sum_k\int_s^t
E \|B_k(r,u((r))\|_{V}^2\,dr
\leq |t-s|M_3
\end{align*}
for any $0\leq s\leq t\leq T$, which proves (ii).  
\end{proof}

\subsection{Approximation spaces 
and operators $\Pi_n$.}                          \label{sec2.2}
Let 
$ V_n \hookrightarrow  H_n\hookrightarrow   V^*_n$
be a normal triple and 
$\Pi_n:  V \to V_n$ be  a bounded linear operator 
for each integer $n\geq0$ 
such that for all $v\in H$ and  
$n\geq0$  
\begin{equation}                                  \label{bpin}
\|\Pi_n v\|_{V_n}
\leq p |v|_V  
\end{equation}
with some constant 
$p$ independent of $v\in V$ and $n$. 
Note that we do not require that  
the maps $\Pi_n$  be orthogonal 
projections on the Hilbert space $H$.

We denote by $\< v,w\>_n$ the duality 
between $v\in V_n$ and $w\in V_n^*$ 
and similarly by
$(h,k)_n$ the inner product of $h,k\in H_n$.  
To lighten the notation, 
let  $\|v\|:=\|v\|_V$ denote the norm 
of $v$ in $V$, $\|v\|_n:=\|v\|_{V_n}$ 
the norm of $v$ in $V_n$, $|u|:=|u|_H$ 
the norm of $u$ in $H$, 
$|u|_n =|u|_{H_n}$ the norm of $u$ in $H_n$,  
and finally  
$|w|_*: = |w|_{V^*}$ and $|y|_{n^*}:=|y|_{V_n^*}$ 
the norm of $w\in V^*$ in $V^*$ and  
the norm of $y$ in $V_n^*$, respectively.

For $r\geq 0$ let $H^r=W^r_2(\mathbb R^d)$ denote the 
closure of $C_0^{\infty}(\mathbb R^d)$ in the 
norm defined by 
$$
|\varphi|_{H^r}^2=\sum_{|\gamma|\leq r}
\int_{{\mathbb R}^d}|D^{\gamma}\varphi(x)|^2\,dx.
$$
In particular,
$H^0=L_2({\mathbb R}^d)$. 
\smallskip

The following basic examples 
 will be used in the sequel.
It describes spaces $V_n, H_n$ and $V^*_n$  and
operators $\Pi_n$ such that condition  
\eqref{bpin} is satisfied.

\begin{Ex}                                          \label{ExWavelets}
{\bf Wavelet approximation.}
{\rm Let $\varphi:\RR\to \RR$ be 
an orthonormal scaling function,  
 i.e., a  real-valued, compactly supported function,  
such that:

(i)   there exists a sequence 
$(h_k)_{k\in \ZZ}\in l^2(\ZZ)$
for which  
$\varphi(x)=\sum_k h_k \varphi(2x-k)$ in $L^2(\RR)$ 
,

(ii) $\int \varphi(x-k) \varphi(x-l) dx 
= \delta_{k,l}$ for any $k,l \in \ZZ$.

\noindent 
We assume that  
the scaling function  
$\varphi$ belongs to the Sobolev space 
$H^s({\mathbb R}):=W^s_2({\mathbb R})$ for sufficiently 
large integer $s>0$. 
\smallskip

For $d>1$, 
$x=(x_1, \cdots, x_d)\in \RR^d$, set 
$\phi(x)=\varphi(x_1)\cdots
\varphi(x_d)$ and for
 $j\geq 0$ and $k\in \ZZ^d$, set 
$\phi_{j,k}(x)=2^{\frac{jd}{2}} \phi(2^jx-k)
\in H^s:=W^s_2({\mathbb R}^d)$.  
For any integer $j\geq 0$, let
$H_j$ denote the closure in 
$L^2(\RR^d)$ of the vector space
generated by $(\phi_{j,k}, k\in \ZZ^d)$ and  
define the 
operator $\Pi_j $ by  
\begin{equation*}                                  \label{Pijwavelet}
\Pi_j f = \sum_{k\in \ZZ^d} 
\big( 
f\, ,\, {\phi}_{j,k}
\big){\phi}_{j,k} \, ,\quad f\in L^2(\mathbb R^d), 
\end{equation*}
$\big(\, ,\, \big)$ denotes the 
 scalar product in $L^2(\RR^d)$.  

Thus we have a sequence 
$H_j\subset H_{j+1}$ of closed subspaces 
of $L^2(\mathbb R^d)$ and orthogonal
projections 
$\Pi_j:L^2(\mathbb R^d)\to H_j$ for $j\geq0$. 
Assume, moreover that $\cup_{j=0}^{\infty}H_j$ 
is dense in $L^2(\mathbb R^d)$ and that 
$\varphi$ is sufficiently regular, such that 
the inequalities 

\begin{align}
(Direct)
\hspace{1cm} & &
\|f-\Pi_j f\|_{H^r}
\leq  C\, 2^{-j(s-r)}\, \|f\|_{H^s}\,, 
\quad \forall f\in H^s, 
                                                 \label{direct}\\
(Converse)
\hspace{1cm}&   & \|\Pi_j f\|_{H^s}\leq  C\,
2^{j(s-r)}\, \|f\|_{H^r},\quad \forall  f\in H^r
                                                  \label{converse}
\end{align} 
holds for fixed integers $0\leq r\leq s$.  
The proof of these inequalities and more information 
on wavelets can be found, e.g., in \cite{C}.  

Fix $r>0$, set $H:=L^2(R^d)$, $V:=H^r=W^r_2(\mathbb R^d)$, 
and identify
$H$ with its dual $H^*$  
by the help of the inner product in $H$.
Then 
$V\hookrightarrow H^*\hookrightarrow V^*$ 
is a normal triple, where 
$H\equiv H^*\hookrightarrow V^*$  is the adjoint of the embedding 
$V\hookrightarrow H$.  
We define $V_n$ as the normed space we get by 
taking the $H^r$ norm on $H_n$.  Since the $H^r$ and $H^0$ 
norms are equivalent on $H_n$, the space $V_n$ is complete,   
and obviously 
$V_n\hookrightarrow H_n\equiv H_n^*\hookrightarrow V_n^*$ 
is a normal triple, where $H_n$ is identified with 
$H_n^*$ via the inner product $(\,,)_n=(\,,)$ in
$H_n$. Note that due to \eqref{converse} we 
have \eqref{bpin} assuming that $\varphi$ is sufficiently smooth. }
\end{Ex}

\begin{Ex}                                       \label{ExFiniteEl} 
{\bf Finite differences -- Discrete
Sobolev spaces.}   
{\rm Consider for fixed
$h\in(0,1)$  the grid 
$$
{\mathbb G}=h{\mathbb Z}^d
=\{(k_1h,k_2h,\dots,k_dh):
k=(k_1,k_2,\dots,k_d)\in\mathbb Z^d\}, 
$$                                                
where $\mathbb Z$ denotes the set of integers. 
Use the notation $\{e_1,e_2,...,e_d\}$ 
for the standard 
basis in $\mathbb R^d$. 
For any integer $m\geq 0$, let $W^{m}_{h,2}$ be  
the set of real valued functions $v$ on 
$\mathbb G$ with 
$$
|v|_{h,m}^2:=\sum_{|\alpha|\leq m}
\sum_{z\in\mathbb G}|\delta^{\alpha}_+v(z)|^2h^d<\infty,
$$ 
where $\delta^0_{\pm i}$ is the identity  and  
$
\delta_{\pm}^{\alpha}
=\delta^{\alpha_1}_{\pm 1}\delta^{\alpha_2}_{\pm 2}
\dots \delta^{\alpha_d}_{\pm d}
$
for multi-indices 
$\alpha=(\alpha_1,\alpha_2\dots,\alpha_d)
\in\{0,1,2,\dots\}^d$ of
length
$|\alpha|:=\alpha_1+\dots+\alpha_d  \geq 1 $ 
is defined for by  
$$
\delta_{\pm i}v(z):=\pm \tfrac{1}{h}(v(z\pm he_i)-v(z)) .
$$ 
We write also $\delta^{\alpha}$ and $\delta_i$ 
in place of $\delta^{\alpha}_{+}$ and $\delta_{+i}$, 
respectively. Then 
$W^m_{h,2}$ with the norm $|\cdot|_{h,m}$ 
is a separable Hilbert space. It is the  
discrete counterpart of the Sobolev space 
$W^{m}_2(\mathbb R^d)$. 
Set $W^{-1}_{h,2}=(W^{1}_{h,2})^{\ast}$, 
the adjoint of $W^{1}_{h,2}$, with its norm 
denoted by $|\cdot|_{h,-1}$.  
It is easy to see  
that 
$
W^{m}_{h,2}\hookrightarrow 
W^{m-1}_{h,2}
$ 
is a dense and continuous embedding, 
\begin{eqnarray}                           \label{15.12.04}              
|v|_{h,m-1}&\leq& |v|_{h,m},               \nonumber \\
|v|_{h,m}&\leq &\frac{\kappa}{h}|v|_{h,m-1} ,   \label{23.06.06.07}
\end{eqnarray} 
for all $v\in W^{m}_{h,2}$, $m\geq0$ 
and $h\in(0,1)$, where $\kappa$ is a constant 
depending only on $d$. Notice that 
for $m\geq1$
$$
\langle v,u\rangle:=\sum_{|\alpha|\leq m}
\sum_{z\in \mathbb G}\delta^{\alpha}v
\delta^{\alpha}u\leq C|v|_{h,m-1}|u|_{h,m+1}
\quad
\text{for all $v,u\in W^{m+1}_{h,2}$}
$$ 
extends to a duality 
product between 
$W^{m-1}_{h,2}$ and $W^{m+1}_{h,2}$, which makes 
it possible to identify $W^{m-1}_{h,2} $ with 
$(W^{m+1}_{h,2})^{\ast}$.
  
Assume that $m>\frac{d}{2}$. Then by Sobolev's theorem 
on embedding 
$W^m_2:=W^m_2(\mathbb R^d)$ into 
${\mathcal C}(\mathbb R^d)$,   
there is a bounded linear operator 
$I:W^m_2(\mathbb R^d)\to {\mathcal C}(\mathbb R^d)$, 
such that $Iu=u$ almost everywhere on $\mathbb R^d$. 
Thus, identifying $u$ with  $Iu$, we can define 
the operator 
$R_h:W^{m}_2(\mathbb R^d)\to W^{m}_{h,2}$ 
by restricting the functions  
$u\in W^{m}_2$  onto 
${\mathbb G}\subset\mathbb R^d$. Moreover,  
due to Sobolev's theorem,  
$$
\sum_{z\in{\mathbb G}}
\sup_{x\in \mathcal I (z)}|u(x)|^2 \, h^d
\leq p^2 |u|^2_{W^m_2}, 
$$
where 
$
\mathcal I (z):=\{x\in{\mathbb R}^d:z_k\leq x_k\leq z_k+h, 
\quad k=1,2\dots d\} 
$
and 
$p$ is a constant depending only 
on $m$ and $d$. Hence obviously 
\begin{equation}                               \label{13.31.03}
|R_h u|_{h,0}^2\leq p |u|^2_{W^m_2}
\quad \text{for all $u\in W^m_2$}.  
\end{equation}                                
Moreover, for every integer $l\geq 0$ 
\begin{equation}                                     \label{2.13}
|R_h u|_{h,l}\leq p |u|_{W^{m+l}_2} 
\quad \text{for all $u\in W^{m+l}_2$}, 
\end{equation}
with a constant $p$ depending only on 
$m$, $l$ and $d$.  Thus setting 
$$
V_n:=W^{m+l}_{h_n,2}, 
\quad 
H_n:=W^{m+l-1}_{h_n,2}, 
\quad 
V_n^{\ast}\equiv W^{m+l-2}_{h_n,2} , 
$$
$$
\Pi_n:=R_{h_n}
$$ 
for any sequence 
$\{h_n\}_{n=0}^{\infty}\subset (0,1)$ 
and any integers $m>\frac{d}{2}$, $l\geq0$ 
we get examples of approximation spaces. 

When approximating differential operators by 
finite differences we need to estimate 
$D_iu-\delta_{\pm i}u$ in discrete Sobolev norms. 
For $d=1$ we can estimate this as follows. 
Let $l\geq0$ be an integer and set 
$z_k:=kh$ for $k\in\mathbb Z$. 
By the mean value theorem there 
exist $z^{\prime}_k$ and $z^{\prime\prime}_k$ 
in $[z_k,z_{k}+lh]$ 
such that $\delta^lDu(z_k)=D^{l+1}u(z^{\prime}_k)$ 
and $\delta^l\delta u(z_k)=D^{l+1}u(z^{\prime\prime}_k)$, 
where $D:=\frac{d}{dx}$.  
Hence 
\begin{eqnarray*}
|\delta^l(Du(z_k)-\delta u(z_k))|^2 &= &
|D^{l+1}u(z^{\prime}_k)-D^{l+1} u(z^{\prime\prime}_k)|^2=
\Big|\int_{z^{\prime}_k}^{z^{\prime\prime}_k}
D^{l+2}u(y)\,dy\Big|^2 \\
&\leq& lh \int_{z_k}^{z_{k}+lh}|D^{l+2}u(y)|^2\,dy 
\end{eqnarray*}
for $u\in C_0^{\infty}(\mathbb R)$. 
Consequently, 
\begin{equation}                                 \label{16.28.03}
|Du-  \delta^l_{\pm}  u|_{h,l}
\leq  l  \, h\, |u|_{W^{l+2}_2(\mathbb R)} 
\end{equation}
for $u\in C^{\infty}_0(\mathbb R)$,  
and hence for all 
$u\in W^{l+2}_2(\mathbb R)$. For $d>1$ by similar 
calculation combined with Sobolev's 
embedding, we get 
that for $m>l+2+\frac{d-1}{2}$
\begin{equation}                                 \label{17.28.03}
|D_iu-\delta_{\pm i}u|_{h,l}
\leq Ch |u|_{W^m_2} 
\end{equation}
for all  $u\in W^m_2$, 
$h\in(0,1)$, where $C$ 
is a constant depending on $l$, 
$m$ and $d$.  }
\end{Ex}

\section{Space discretization}
\subsection{Description  of the scheme} 

Consider  for each 
integer $n\geq1$ the problem 
\begin{equation}                                  \label{11.03.04} 
du^n(t)=A^n(t,u^n(t))\,dt
+\sum_kB^n(t,u^n(t))\,dW^k(t), \quad u^n(0)=u^n_0, 
\end{equation}
in a normal triple 
$ V_n \hookrightarrow  H_n\hookrightarrow   V^*_n$, 
satisfying
the conditions  of section \ref{sec2.2}, 
where $u^n_0$ is an $H_n$-valued 
$\mathcal F_0$-measurable 
random variable, and $A^n$ and $B^n=(B^n_k)$ 
are $\mathcal P\otimes{\mathcal B}(V_n)$-measurable 
mappings from $[0,\infty)\times \Omega\times V_n$ 
into $V_n^*$ and $H_n^{d_1}$, respectively.

\begin{assumption}                            \label{assumption S}
The operators $A^n$ and $B^n$ satisfy the following 
conditions. 

\noindent
{\bf (S1)}  
 ({\it Strong monotonicity})  
{\it There exist constants $\lambda>0$ and $L$ such that 
for all $n\geq1$ 
almost surely 
\begin{align*}  2\< u-v , A^n(t,u)-A^n(t,v)\>_n
 &+\sum_k|B^n_k(t,u)-B^n_k(t,v)|_{H_n}^2
+ \lambda  \|u-v\|_{V_n}^2            \\
 &\leq L |u-v|_{H_n}^2
\quad
\text{for all $t\in [0,T]$, $u,v\in V_n$}.                         
\end{align*}
}
{\bf (S2)} (Growth condition)
{\it   Almost surely 
\begin{equation*}                                          
\label{boundAn}
 |A^n(t,v)|_{V_n^*}^2 \leq K_{1} \, 
\|v\|_{V_n}^2 + f^n(t), 
\quad
|B^n(t,v)|_{H_n}^2 \leq K_{2} \, 
\|v\|_{V_n}^2 + g^n(t)\,                                            
\end{equation*}
for all $t\in[0,T]$, $v\in V_n$ and $n\geq1$, 
where $K_1$, $K_2$ are  constants, 
independent of $n$, and  
$f^n$ 
and $g^n$   are  non-negative 
stochastic processes such that 
\begin{equation*}                                   \label{bfn}
\sup_n E\int_0^T f^n(t)\, dt =:M_1 <\infty\; , 
\quad
\sup_n E\int_0^T g^n(t)\, dt =:M_2 <\infty.
\end{equation*}

\noindent 
{\bf (S3)} (Hemicontinuity of $A^n$) 
For every $n\geq 1$, the operators $A^n$ 
are hemicontinuous in
$v\in V_n$,  i.e., almost surely 
$$
\lim_{\varepsilon\to0}
\langle A^n(t,v+\varepsilon u)\,,w\rangle_n=
\langle A^n(t,v)\,,w\rangle_n
$$ 
for all $t\in[0,T]$, $v,u,w\in V_n$. 

\noindent
{\bf (S4)} (Lipschitz condition on $B^n$)  
Almost surely 
\begin{equation*}                                    \label{lipBn}
\sum_k |B^n_k(t,u) - B^n_k(t,v)|_{H_n}^2  
\leq L_B\, \|u-v\|_{V_n}^2
\end{equation*} 
for all  $t\in [0,T]$ and $u,v\in V_n$. }
\smallskip

\end{assumption}

 The solution to \eqref{11.03.04} is understood 
in the sense of Definition \ref{sol}.
Notice that  {\bf (S1)} - {\bf (S2)} 
imply the  coercivity condition 
$$
2\<v , A^n(t,v)\>_n
+\sum_k|B^n_k(t,v)|_{H_n}^2
+ \tfrac{\lambda}{2}  \|v\|_{V_n}^2
\leq  C\, \big(|v|_{H_n}^2+f^n(t)+g^n(t)\big)
$$
with a constant $C$ depending on $\lambda$,  
$L$ and $K_2$. 

Thus by Theorem \ref{existu} 
the conditions {\bf(S1)}--{\bf(S3)} 
ensure 
the existence of a unique solution $u^n$ to 
\eqref{11.03.04}, and if 
\begin{equation}                                \label{In}
\sup_n E|u_0^n|_{H_n}^2  <\infty, 
\end{equation}
then
\begin{align}                                   \label{bun}         
E\sup_{0\leq t\leq T} & |u^n(t)|_{H_n}^2 
+  E\int_0^T \|u^n(t)\|_{V_n}^2\, dt \nonumber \\
&\leq C \,  \sup_n \left(E|u_0^n|^2_{H_n}
+E\int_0^T\big( f^n(t)+g^n(t)\big) \,dt\right)<\infty, 
\end{align}
where $C$ is a constant depending only on 
$\lambda$, $L$ and $K_2$.   

\subsection{Rate of convergence of the scheme}

We want to approximate $\Pi_n u$ by 
$u^n$. In order to estimate the accuracy
of this approximation  
we need to relate the operators
$A$ and $B$ to $A^n$ and  $B^n$, respectively. 
Therefore we assume the regularity condition 
{\bf (R3)} from Assumption \ref{assumption R} 
and make the following {\it consistency} assumption.

\smallskip

\noindent 
{\bf Condition (Cn)} {\it (Consistency)}  
{\it There exist 
a sequence $(\varepsilon_n)_{n\geq1}$ 
of positive  
numbers and a sequence  $(\xi^n)_{n\geq1}$  
of non-negative adapted processes such
that
$$
 \sup_n E\int_0^T \xi^n(t)\, dt\leq M <+\infty ,
$$
 and almost surely  
$(t,\omega)\in[0,T]\times\Omega$
\begin{align*}                                              \label{consisn}
|\Pi_n A(t,v) -  A^n\big(t,  \Pi_n v\big)|_{V_n^*}^2 
+  & \sum_k 
|\Pi_n B_k(t,v) -B^n_k\big(t,\Pi_n v\big)|_{H_n}^2 \nonumber \\
&  \leq  \varepsilon_n^2 
\big(|v|_{\mathcal V}^2 + \xi^n(t) \big) 
\end{align*}
for all $t\in[0,T]$ and $v\in\mathcal V$.}

\begin{Th}                                 \label{raten}
Let  Assumption \ref{assumption S},  
the regularity conditions 
{\bf(R1)} and {\bf(R3)} from 
Assumption \ref{assumption R}, and the 
consistency condition 
{\bf (Cn)}   
hold. Assume furthermore 
$\sup_n E|u_0^n|_{H_n}^2 <+\infty$. 
Then for $e^n(t):=\Pi_n u(t)-u^n(t)$,  
\begin{equation}                            \label{nEsup}
E \sup_{0\leq t\leq T} 
|e^n(t)|_{H_n}^2 
+ E\int_0^T\|e^n(t)\|_{V_n}^2
dt 
 \leq C_1 E|e^n(0)|_{H_n}^2
+C_2(r_1 + M)\varepsilon_n^2  
\end{equation}
holds for all $n\geq1$, 
where $C_1=C_1(\lambda, L,T)$  and 
$C_2=C_2(\lambda,L,L_B,T)$ are constants.  
\end{Th}
\begin{proof}  From equation \eqref{u} 
we deduce that for every $n\geq 1$,
$$
\Pi_n u(t)=\Pi_n u_0 
+ \int_0^t \Pi_n A\big(s,u(s)\big)\, ds +
\sum_k \int_0^t\Pi_n B_k\big(s,u(s)\big) dW^k(s).
$$
Using It\^o's formula 
\begin{equation}                              \label{Iton}
|e^n(t)|_n^2 
= |e^n(0)|_n^2 + \sum_{i\leq 3}
I_i(t),
\end{equation}
where
\begin{eqnarray*}
I_1(t)&=&2\int_0^t 
\langle e^n(s)\,,\Pi_n A
(s,u(s)\big) -A^n\big(s, u^n(s))
\rangle_n ds \, ,                                                \\
I_2(t)&=&2\sum_k \int_0^t \big(e^n(s)\,, 
\Pi_n B_k \big(s, u(s)\big) -
B^n_k\big( s,u^n(s)\big) \big)_n dW^k(s),                         \\
I_3(t)&=& \sum_k \int_0^t \big| 
\Pi_n B_k \big(s, u(s)\big) -
B^n_k\big( s,u^n(s)\big) \big|_n^2  ds\, .
\end{eqnarray*}
We first prove
\begin{equation}                                              \label{nsupE}
\sup_{0\leq t\leq T} E  |e^n(t)|_{n}^2  
+ E \int_0^T\!\! \|e^n(t)\|_n^2
dt  
\leq  C_1E|e^n(0)|_n^2  
 +   C_2(r_1 + M)\varepsilon_n^2,   
\end{equation}
where 
 $C_1=C_1(\lambda,L,T)$ and $C_2=C_2(\lambda, L,
L_B, T)$  are constants.   
The strong monotonicity condition  {\bf (S1)} 
from Assumption \ref{assumption S}  
implies 
\begin{equation}                                               \label{I13n}
 I_1(t)+I_3(t)\leq - \lambda 
\int_0^t \|e^n(s)\|_n^2\, ds 
+
L\int_0^t |e^n(s)|_n^2\, ds + \sum_{i=1,2}
R_i(t),
\end{equation}
where
\begin{eqnarray*}
R_1(t)&=&
\int_0^t 2\langle e^n(s)\, ,\, 
\Pi_n A(s,u(s)) -A^n(s, \Pi_n u(s))\rangle_n
\, ds, \\
R_2(t) &= & \sum_k \int_0^t 
\Big[
| \Pi_n B_k(s,u(s)) -B^n_k(s,u^n(s))|_n^2  \\
&& - |B^n_k(s,\Pi_n u(s)) - B^n_k(s,u^n(s))|_n^2 
\Big]\, ds.
\end{eqnarray*}
Schwarz's inequality 
and  the consistency condition {\bf (Cn) }  
imply that  for every $n\geq 1$ 
and $t\in [0,T]$, 
\begin{align}                                      \label{R1n}
|R_1(t)|& \leq  \tfrac{\lambda}{3} 
\int_0^t \|e^n(s)\|_n^2\, ds 
+ \tfrac{3}{\lambda}
\int_0^t |\Pi_n A\big( s,u(s)\big) - A^n
\big(s, \Pi_n u(s)\big) |_{n*}^2\, ds               \nonumber \\
& \leq \tfrac{\lambda}{3} \, 
\int_0^t \|e^n(s)\|_n^2\, ds 
+ \tfrac{3}{\lambda}\, 
\varepsilon^2_n  \,  
\int_0^t \big(|u(s)|_{\mathcal V}^2  
+ \xi^n(s) \big)\, ds .
\end{align}
Schwarz's inequality, 
the consistency condition {\bf (Cn)}, 
and the Lipschitz condition 
{\bf (S4)} from Assumption \ref{assumption S}  
 yield that for every $\alpha >0$,
\begin{align}
|R_2&(t)|=  \sum_k \int_0^t 
\Big[ 
|\Pi_n B_k
\big(s,u(s)\big) - B_k^n 
\big(s, \Pi_n u(s)\big)|_n^2
                                                          \nonumber\\
&\quad 
+ 2\Big( 
\Pi_n B_k\big(s, u(s)\big) - B^n_k\big(s, \Pi_n u(s)\big) 
\, ,\, B ^n_k\big( s,\Pi_n u(s)\big) -
B^n_k\big(s, u^n(s)\big) 
\Big)_n
\Big]
\, ds                                                     \nonumber\\
& \leq ( 1+\tfrac{1}{\alpha}) 
\int_0^t \sum_k \big|\Pi_n B_k\big(s, u(s)\big) -
 B^n_k\big(s,\Pi_n u(s)\big) \big|_n^2 ds                 \nonumber\\
&\qquad 
+ \alpha 
\int_0^t \sum_k \big| B^n_k\big(s, \Pi_n u(s)\big)
- B^n_k \big(s, u^n(s)\big)\big|_n^2 \, ds                \nonumber\\
& \leq ( 1+\tfrac{1}{\alpha}) \, \varepsilon_n^2 
\int_0^t \big(|u(s)|_{\mathcal V}^2 + \xi^n(s)\big)\, ds  
+ \alpha L_B 
\int_0^t \|e^n(s)\|_n^2\, ds.                             \label{R2n}
\end{align}
Thus, for 
$\alpha L_B\leq \frac{\lambda}{3}$, 
taking expectations in \eqref{Iton} and
\eqref{I13n}-\eqref{R2n} 
and using {\bf (S1) } again, we deduce that
$$
E|e^n(t)|_n^2 
+\tfrac{\lambda}{3} 
E\int_0^t \|e^n(s)\|_n^2\, ds              
\leq
L E\int_0^t |e^n(s)|_n^2 ds 
+E|e^n(0)|_n^2 
+C(r_1 + M) \varepsilon_n^2,
$$
where $C=C(\lambda,L_B)$ is a constant. 
Since by
\eqref{bounduH}  and \eqref{bun} 
$$
\sup_{0\leq t\leq T} E|e^n(t)|_n^2<+\infty,
$$
Gronwall's lemma  gives 
$$
\sup_{0\leq t\leq T}  
E|e^n(t)|_n^2
\leq 
e^{LT}\left(C(r_1 + M)\varepsilon_n^2
+ E|e^n(0)|_n^2
\right), 
$$ 
which in turn  yields \eqref{nsupE}. 
\smallskip
We now prove \eqref{nEsup}.  
From \eqref{nsupE}--\eqref{R2n} we deduce 
\begin{align}                                   \label{In24}     
E\sup_{0\leq t \leq T}\big(I_1(t)+I_3(t)\big)
& \leq LE\int_0^T \!\!\!
|e^n(s)|_n^2\, ds 
+\tfrac{2\lambda}{3} 
E \int_0^T  \!\!\! \|e^n(s)\|_n^2 ds 
+C_2 (r_1+M) \, \varepsilon_n^2              
                                                 \nonumber\\
& \leq  C_1E|e^n(0)|_n^2+
C_2(r_1 + M)\varepsilon_n^2 .                  
\end{align}
(Notice that by taking the supremum 
in both sides of \eqref{I13n} we cannot make use of the 
term with coefficient  $-\lambda$ 
in the right-hand side of \eqref{I13n}. This is why 
$2\lambda/3$ appears here  as the sum of 
$\lambda/3$ from \eqref{R1n} and $\alpha L_B\leq\lambda/3$ 
from \eqref{R2n}.) By Davies'
inequality,   
\eqref{bounduH},  
the Lipschitz condition  {\bf (S4)} on $B^n$, 
the consistency condition  
{\bf (Cn)} and by the strong monotonicity 
condition {\bf (S1)}, 
\begin{align}
& E\sup_{0\leq t\leq T}  |I_2(t)| 
\leq 6E\Big( \int_0^T 
\sum_k
\Big| \Big(e^n \, ,\Pi_n B_k(u)
-B^n_k(u^n)\Big)_n\Big|^2ds\Big)^{\frac{1}{2}} \nonumber \\ 
& \; \leq 
6E\Big\{ \sup_{0\leq t\leq T}
 |e^n(t)|_n                                  
\Big( 
\int_0^T\sum_k |\Pi_n B_k( u) -
B^n_k\big(u^n\big)|_n^2\,ds \Big)^{\frac{1}{2}} \Big\}
\nonumber \\ & \;  \leq \frac{1}{2} 
E\sup_{0\leq t\leq T} 
|e^n(t)|_n^2                                              
                                                      \nonumber\\  
& \quad +   36  E\sum_k \! \int_0^T \! \! 
\Big[|\Pi_n B_k\big(u\big) - B_k^n\big(\Pi_n u\big)|_n^2
+ |B^n_k\big(\Pi_n u\big) -B^n_k(u^n)|_n^2 \Big]\, ds    
\nonumber \\
                                        \label{In3}
&\;  \leq \tfrac{1}{2} 
E\sup_{0\leq t\leq T} 
|e^n(t)|_n^2 
+   36  L_B E\int_0^T \|e^n(s)\|_n^2 ds      
+  36 (r_1 + M)\varepsilon_n^2,  
\end{align}
where the argument $s$ is omitted from most 
integrands. 
 Thus,  relations  \eqref{Iton}, \eqref{In24}, 
\eqref{In3} and  \eqref{nsupE}
 yield
$$
\frac{1}{2} E\sup_{0\leq t\leq T} |e^n(t)|_n^2 
\leq C_1E|e^n(0)|_n^2 +C_2(r_1 + M)\varepsilon_n^2 , 
$$
with some 
constants $C_1=C_1( L,T)$   and 
$C_2=C_2(\lambda, L,L_B, T)$,  
which completes the proof of 
\eqref{nEsup}.
\end{proof}

\subsection{Example}                           \label{ExSpace}
Consider the normal triples 
$$
V\hookrightarrow H^*\hookrightarrow V^*, 
\quad V_n\hookrightarrow H_n^*\hookrightarrow V_n^*
$$
with the orthogonal projection 
$\Pi_n: H=L^2(\mathbb R^d)\to H_n$  
from 
Example \ref{ExWavelets}, where 
$V=W_2^r(\mathbb R^d)$ with $r>0$. 
Set $\mathcal H=W^{r+\rho}_2(\mathbb R^d)$ and 
$\mathcal V=W^{r+l}_2(\mathbb R^d)$ for some 
$l>\rho\geq0$. 

Let  $A$ and   $B=(B_k)$ be 
$\mathcal P\otimes\mathcal B(V)$-measurable mappings 
from $[0,\infty[\times \Omega\times V$ 
into $V^*$ and $H^{d_1}$, respectively, 
satisfying Assumptions \ref{E} and 
\ref{assumption R}. 
 For $(t,\omega)\in[0,T]\times\Omega$ 
let $A^n(t,\omega,\cdot):V^n\to V_n^*$ and 
$B^n(t,\omega,\cdot):V^n\to H_n^{d_1}$ be defined 
by 
\begin{equation}                             \label{exoperatorspace}
\<A^n(t,u,\omega)\,,v\>_n=\<A(t,u,\omega)\, ,v\> 
\quad \text {and $B^n_k(t,\omega,u)=\Pi_n B_k(t,\omega,u)$}
\end{equation}
for all $u,v\in V_n$, where $\<\,,\>_n$ denotes 
the duality between $V_n$ and $V_n^*$. Then 
 it is easy to see that due to conditions 
{\bf(1)}, {\bf(2)} and {\bf(3)} in Assumption \ref{E}, 
the operators 
$A^n$ and $B^n$ satisfy  
{\bf (S1)}, {\bf (S2)} and 
{\bf (S3)} in Assumption \ref{assumption S}, respectively. 
Furthermore, taking into account 
Remark \ref{remark 29.09.08} it is obvious that 
{\bf(S4)} holds. 
Assume the regularity condition {\bf (R3)} from 
Assumption \ref{assumption R}. Then by virtue of 
the definition of $\Pi_n$, $A^n$ and $B^n$,  
due to 
Lipschitz conditions  {\bf(3)}  in Assumption \ref{E} and
\eqref{29.09.08}  in Remark \ref{remark 29.09.08}, 
we have, recalling the direct inequality  
\eqref{direct},  
\begin{equation*}
|\Pi_n A(t,u) -A^n(t,\Pi_n u)|_{V_n^*}^2 
+ \sum_k |\Pi_n B_k(t,u) - B^n_k(t, \Pi_n u)|_{H_n}^2
\end{equation*}
\begin{equation*} 
\leq |A(t,u)-A^n(t,\Pi_n u)|_{V*}^2 
+ \sum_k |B_k(t,u) - B^n_k(t, \Pi_n u)|_{H}^2
\end{equation*}
\begin{equation*} 
\leq C(L_1 + L_2)\, 2^{-2nl}\, 
|u|_{\mathcal V}^2 
\end{equation*}
almost surely for all $t\in[0,T]$ and $u\in {\mathcal V}$, 
which  
yields {\bf (Cn)} with $\xi^n:=0$ and 
$\varepsilon_n:=C(L_1+L_2)2^{-nl}$. 
In the last section we will give 
examples of operators such that 
 Assumption \ref{assumption R} holds.

\section{Implicit space-time discretizations}
\subsection{Description  of the scheme}

For a fixed integer $m\geq1$ set $\tau:=T/m$ 
and $t_i=i\tau$ for $i=0, \cdots, m$.
Let $V_n\hookrightarrow H_n \hookrightarrow V_n^*$
 satisfy the conditions in section \ref{sec2.2}. 
 Given a  $V_n$-valued 
 $\mathcal F_0$-measurable random variable  
$u_0^{n,\tau}$ and  
${\mathcal F}_{t_i}\otimes {\mathcal B}(V_n)$-measurable  
mappings 
$$
A^{n,\tau}_{j} : \Omega \times V_n  \to V_n^*\;   
\quad \mbox{\rm and} 
\quad  B_{k,i}^{n,\tau} : \Omega\times V_n 
\to H_n,\;  \mbox{\rm for }\,  k=1, \cdots, d_1,    
$$  
$j=1,\dots,m$ and $i=0,\dots m-1$, 
consider for each $n$   the
system  of equations 
 \begin{equation}                                       \label{unspi}
u^{n,\tau}_{i+1}= u^{n,\tau}_i
+ \tau\,  A^{n,\tau}_{i+1}
\big(u^{n,\tau}_{i+1}\big)                            
+\sum_k {B}^{n,\tau}_{k,i} 
\big(u^{n,\tau}_i\big)\,
\big(W^k(t_{i+1})- W^k({t_i})\big), 
\end{equation}
$i=0,\dots,m-1$, for $V_n$-valued  
$\mathcal F_{t_{i}}$-measurable  
random variables  $u^{n,\tau}_{i}$, 
$i=1,\dots,m$. 
 
\smallskip
\begin{assumption}                           \label{assumption ST} 
For almost all $\omega\in\Omega$ the operators 
$A^{n,\tau}_{j}$ and 
$B^{n,\tau}_{k,i}$ satisfy the following 
conditions for all $j=1,\dots,m$, $i=0,\dots,m-1$,

\noindent
{\bf (ST1)} (Strong monotonicity) 
There exist constants $\lambda>0$ and $L\geq0$ such that 
\begin{align}                                           
 2\big\langle u-v\, ,\,  
A^{n,\tau}_j(u)-A^{n,\tau}_j(v)\big\rangle_n &
 +\sum_k 
\big| 
B^{n,\tau}_{k,j}(u)- B^{n,\tau}_{k,j}(v)
\big|_{H_n}^2
                                                     \nonumber \\
& \leq  - \lambda \|u-v\|_{V_n}^2
 + L\,
|u-v|_{H_n}^2                                        \label{monotonspi}
\end{align}
for all $u,v\in V_n$, $m\geq1$, $n\geq0$.
  
\noindent  
{\bf (ST2) } (Growth condition on  $A^{n,\tau}_i$ and 
$B^{n,\tau}_i$)   
There is a constant $K$ such that  
\begin{equation*}                                 \label{growthAspi} 
\big|A^{n,\tau}_{j}(u)\big|_{V_n^*}^2
\leq K\|u\|_{V_n}^2+f^{n,\tau}_{j}, 
\quad 
 \sum_k  \big|B^{n,\tau}_{k,i}(u)\big|_{H_n}^2
\leq K\|u\|_{V_n}^2+g^{n,\tau}_{i}  
\end{equation*}
for all $u\in V_n$, $m\geq1$, $n\geq0$, 
 where 
 $f^{n,\tau}_{j}$ and $g^{n,\tau}_{i}$ are 
non-negative
random variables, such that 
\begin{equation*}                                       \label{fnspi}
\sup_{n,m} 
\sum_{j} \tau Ef^{n,\tau}_{i}\leq M <+\infty, 
\quad 
\sup_{n,m}
\max_{i}Eg^{n,\tau}_{i} \leq M<+\infty.
\end{equation*}

\noindent  
{\bf (ST3)}  
(Lipschitz condition on $A^{n,\tau}_j$ )  
There exists  a constants $L_1$   
such that 
\begin{equation}                                      \label{LipAAspi}
\big|
A^{n,\tau}_{j}(u)- A^{n,\tau}_{j}(v)
\big|_{V_n^*}^2
\leq L_1 \|u-v\|_{V_n}^2          
\end{equation} 
for all $u,v\in V_n$, $m\geq1$, $n\geq0$.   
\end{assumption}

\begin{rk}                                         \label{remark Lip}
Clearly, conditions {\bf(ST1)} and {\bf(ST3)} imply 
the Lipschitz continuity of $B^{n,\tau}_{k,i}$ in $v\in V_n$, 
i.e., there is a constant $L_2=L_2(L,\lambda,L_1)$ 
such that almost surely 
\begin{equation}                                     \label{lipBnm}
\sum_k \big|                       
B^{n,\tau}_{k,j}(u)- B^{n,\tau}_{k,j}(v)\big|_{H_n}^2
\leq  L_2\, \|u-v\|_{V_n}^2
\end{equation}
for all $u,v\in V_n$, $n\geq 1$, $m\geq 0$ 
and  $j=1,\cdots,m$. 
\end{rk}

\begin{rk}                                         \label{rk41}
Conditions   {\bf (ST1)}--{\bf (ST2)}   
imply that 
almost surely
\begin{equation*}                                  \label{coercivAntau}
 2\< u,  A^{n,\tau}_{j}(u)\>_n  
+ \sum_k |B^{n,\tau}_{k,j}(u)|_{H_n}^2 
\leq
-\tfrac{\lambda}{2}\, \|u\|_{V_n}^2 
+ L|u|_{H_n}^2 + C(f^{n,\tau}_{j}+g^{n,\tau}_{j}) 
\end{equation*}
for all $u\in V_n$, $n\geq0$, $m\geq1$ and $j=1,\dots,m$, 
where  $C=C(\lambda,K)$ is a constant. The Lipschitz 
condition {\bf (ST3)}  
obviously 
implies that $A^{n,\tau}_{i}$ is hemicontinuous. 
\end{rk}

\begin{proposition}                          \label{22.11.04}
Let Assumption \ref{assumption ST} hold. 
Assume $E\|u_0^{n,\tau}\|^2_{V_n}<\infty$ 
for all $n\geq0$ and $m\geq 1$.   
Then 
for $\tau<1/L$  
equation \eqref{unspi} has a unique $V_n$-valued solution 
$(u^{n,\tau}_{i})_{j=1}^m$, such that $u^{n,\tau}_{j}$ 
is $\mathcal F_{t_j}$-measurable and  
$E\|u^{n,\tau}_{j}\|^2_{V_n}$ 
is finite for each $j$, $n$. (Here $1/L:=\infty$ 
if $L=0$.)
\end{proposition}
\begin{proof} Equation \eqref{unspi} 
can be rewritten as 
\begin{equation}                                \label{10.04.04}
D_{i+1} (u^{n,\tau}_{i+1})= u^{n,\tau}_i                            
+\sum_k {B}^{n,\tau}_{k,i} 
\big(u^{n,\tau}_i\big)\,
\big(W^k(t_{i+1})- W^k({t_i})\big), 
\end{equation}
where $D_i:V_n\to V_n^*$ is defined by 
$D_i(v)=v-\tau A^{n,\tau}_{i}(v)$ for each 
$i=1,2,\cdots m$. Due to  
Assumption \ref{assumption ST}
and Remark \ref{rk41}
the operator $D_i$ satisfies the 
assumptions (monotonicity, coercivity, linear growth 
and hemicontinuity) of Proposition 3.4  in
\cite{GyMi}  with $p=2$. By virtue 
of this proposition, for $\tau<1/L$,  equation 
\eqref{10.04.04} has a unique $V_n$-valued 
$\mathcal F_{t_{i+1}}$-measurable solution 
$u^{n,\tau}_{i+1}$ 
for every
given 
$V$-valued $\mathcal F_{t_{i}}$-measurable 
random variable $u^{n,\tau}_{i}$,  and 
\begin{eqnarray*}
E\|u^{n,\tau}_{i+1}\|^2_{V_n}
&\leq  &C\, E
\Big( 1+ f_i^{n,\tau}+g_i^{n,\tau}   
+\Big|\sum_k {B}^{n,\tau}_{k,i} 
(u^{n,\tau}_i)\,
(W^k(t_{i+1})- W^k({t_i}))\Big|^2\Big) \\
&\leq & C
\Big( 1+  Ef_i^{n,\tau} +Eg_i^{n,\tau} 
+\sum_k \tau  E|{B}^{n,\tau}_{k,i} 
(u^{n,\tau}_i)|^2_n \Big)    \\
&  \leq & C\big(1+ Ef_i^{n,\tau} +Eg_i^{n,\tau} 
+K \tau E\|u^{n,\tau}_{i}\|^2_{V_n}
+ \tau\, Eg_i^{n,\tau}\big), 
\end{eqnarray*}
where $C=C(\lambda,\tau)$ is a constant. 
Hence induction on $i$ concludes the proof. 
\end{proof}

\subsection{Rate of convergence of the implicit  scheme} 

Let Assumption \ref{assumption R} on the regularity of 
equation \ref{v} and its solution $u$ hold. 
We relate 
the operators $A(t_i,.)$ and
$A^{n,\tau}_i$ as well as the operators 
$B_k(t_i,.)$ and $B^{n,\tau}_{k,i}$     
by the following consistency assumption. 
\smallskip

\noindent   
{\bf Condition  (Cn$\tau$)} {\it (Consistency)}   
{\it There exist constants $\nu\in ]0,\frac{1}{2}]$, 
$c\geq0$,  
a sequence of numbers $\varepsilon_n\to 0$,
such that
almost surely 
\begin{align}
|\Pi_n A(t_j,u)-A^{n,\tau}_{j}(\Pi_n u)|_{V_n^*}^2  
&\leq c\big(|u|_{\mathcal V}^2 
+ \xi_{j}^{n,\tau}\big) \big(\tau^{2\nu}  +
\varepsilon_n^2\big)   ,
                                                          \nonumber \\
\sum_k
\left|
\Pi_nB_k(t_i,u) -B^{n,\tau}_{k,i}(\Pi_n u)
\right|_{H_n}^2
 &\leq
c (|u|_{\mathcal V}^2 
 + \eta_i^{n,\tau})\big(\tau^{2\nu}  + \varepsilon_n^2\big)   
                                                           \nonumber
\end{align}
for all $j=1,\dots m$, $i=0, \cdots, m-1$ and 
$u\in {\mathcal V}$, where $\xi_j^{n,\tau}$ and 
$\eta_i^{n,\tau}$ are non-negative random 
variables such that  
$$
\sup_{n,m}\sum_{j}   \tau E  \xi^{n,\tau}_j\leq M, 
\quad
\sup_{n,m}\sum_{i}  \tau E \eta^{n,\tau}_i\leq M. 
$$
}

\begin{Th}                                                 \label{th1spi}
 Let 
Assumptions \ref{assumption R} 
and \ref{assumption ST} 
as well as  condition {\bf(Cn$\tau$)} hold. Assume 
\begin{equation}                                         \label{Inm}
\sup_{n,m}E\|u_0^{n,\tau}\|^2_{V_n}\leq M. 
\end{equation}
Set $e_i^{n,\tau}=\Pi_nu(t_i)-u^{n,\tau}_i$. 
 Then for $\tau<1/L$ and 
$n\geq 0$ 
\begin{eqnarray}                                             
&& E\max_{0\leq i\leq m}   
|e^{n,\tau}_i|_{H_n}^2 +
\sum_{1\leq i\leq m}\!  
\tau \, E\|e^{n,\tau}_i\|_{V_n}^2 \nonumber \\
&& \qquad\qquad                                      \label{ratentau}    
\leq C_1  E|e_0^{n,\tau}|_{H_n}^2   
+   C_2( \tau^{2\nu}+ \varepsilon_n^2)
(r_1+r_2+M), 
\end{eqnarray} 
where $C_1=C_1(\lambda, L,T)$ and 
$C_2=C_2(\lambda, L,K,T,p,c,L_1,L_2)$ are constants.
\end{Th}

\begin{proof} We fix $n$, $\tau$, and 
to ease notation we write $e_i$, $A_i$ and $B_{k,i}$ 
in place of $e_i^{n,\tau}$, $A_i^{n,\tau}$ and 
$B_{k,i}^{n,\tau}$, respectively. Similarly, 
we  often 
use $u_i$ in place of $u^{n,\tau}_i$ for $i=1,2,\cdots,m$. 
Then for any $i=0,\cdots, m-1$, 
\begin{align*}
|e_{i+1}|_n^2 &  -|e_i|_n^2                      
=2\, \int_{t_i}^{t_{i+1}}
\big\langle e_{i+1}\, ,
 \Pi_n A(s, u(s)) -  A_{i+1}(u_{i+1})
\big\rangle_n\, ds                                        \nonumber \\
&\quad 
+ 2  \sum_k
\int_{t_i}^{t_{i+1}}\big(e_{i+1}\,, F_k(s)\big)_n\,
dW^k(s)
                                                             \nonumber\\
&\quad  
- \Big| 
\int_{t_i}^{t_{i+1}}
\big[
\Pi_n A(s,u(s)) - A_{i+1}(u_{i+1})
\big]\,ds                                                           
+
\sum_k  \int_{t_i}^{t_{i+1}} F_k(s) \, dW^k(s)
\Big|^2_n                                                     \nonumber\\
= & 2\, \int_{t_i}^{t_{i+1}}
\big\langle e_{i+1}\, ,
 \Pi_n A(s, u(s)) -  A_{i+1}( u_{i+1})
\big\rangle_n\, ds                                            
                                                             \nonumber\\  
&\quad +
\Big| 
\sum_k \int_{t_i}^{t_{i+1}}F_k(s)\,dW^k(s)
\Big|^2_n 
+ 2\sum_k \int_{t_i}^{t_{i+1}} 
\big(e_i\, ,\,F_k(s)\big)_n\,dW^k(s)
\\
&\quad 
-\Big| \int_{t_i}^{t_{i+1}}
\big[ \Pi_nA(s,u(s)) - A_{i+1}(u_{i+1})
\big]\, ds  \Big|_n^2, 
\end{align*}
where for $k=1, \cdots, d_1$ one sets  
\begin{equation*}                                          \label{Fk}
 F_k(s)
=\Pi_n B_k\big(s,u(s)\big) - {B}_{k,i}(  u_i^{n,\tau}), 
\quad s\in]t_i,t_{i+1}], \quad i=0,1,\cdots,m-1 .
\end{equation*}
Summing up for $i=0, \cdots,\,  l-1$, we obtain
\begin{equation}                                           \label{diffspi}
|e_l|_n^2
\leq |e_0|_n^2 
+  2
\sum_{0\leq i<l}\int_{t_i}^{t_{i+1}} \!\!
\langle e_{i+1}\, ,\Pi_n A(s, u(s)) - A_{i+1}(u_{i+1})
\rangle_n \, ds  +Q(t_l)+I(t_l), 
\end{equation}
where 
\begin{align*}                                            
Q(t_l)= &\sum_{0\leq i<l}\Big| 
\sum_k \int_{t_i}^{t_{i+1}}F_k(s)\,dW^k(s)
\Big|^2_n , \\ 
I(t_l)= & 2\sum_k \int_{0}^{t_l} 
\big(e(s)\, ,\,F_k(s)\big)_n\,dW^k(s), 
\; e(s):=e_i  
\;  \text{for $s\in]t_i,t_{i+1}]$, 
$i= 0, \cdots, m$}.                               
\end{align*}
First we show 
\begin{equation}                                      \label{rate1ntau}
\sup_{0\leq l\leq m}E|e_l|_n^2 +
E \sum_{1\leq i\leq m} \tau 
\|e_i\|_n^2 
\leq C_1  E|e_0|_n^2 +  
C_2( \tau^{2\nu}+ \varepsilon_n^2)
(r_1+r_2+M), 
\end{equation} 
where $C_1=C_1(\lambda, L,T)$ and 
$C_2=C_2(\lambda, L,K,T,p,c,L_1,L_2)$ are constants. 
To this end we take expectation in both sides of 
\eqref{diffspi} and use  
the strong monotonicity condition 
{\bf (ST1)}  
from Assumption \ref{assumption ST} 
to get  
\begin{align}                                                     \nonumber
 E|e_l|^2_n &\leq E|e_0|_n^2                                                               
+ 2E\sum_{0\leq i<l}  \tau 
\big\langle e_{i+1}\, ,
A_{i+1} (\Pi_n u(t_{i+1})) - A_{i+1}
(u_{i+1})\big\rangle_n
                                                              \nonumber \\
& \quad
+ E\sum_{0\leq i < l-1}\sum_k  \tau\, 
\big|  B_{k,i+1}(\Pi_n u(t_{i+1}))  
- B_{k,i+1}
(u_{i+1})\big|_n^2
+\sum_{1\leq j\leq 3} S_j                                       \nonumber \\
& \leq  
E|e_0|_n^2   
- \lambda \sum_{1\leq i \leq l} \tau  E\|e_{i}\|_n^2                                            \
+L\sum_{1\leq i\leq l} \tau   E|e_{i}|_n^2
+\sum_{1\leq j\leq 3} S_j                                     \label{majoE1spi}
\end{align}
for $l=1,\cdots,m$, where
\begin{eqnarray*}
S_1&=&2\, \sum_{1\leq i\leq l}
E\int_{t_{i-1}}^{t_i}
\langle e_i\, ,\, \Pi_n
A(s,u(s))-A_i(\Pi_nu(t_i))\rangle_n\,ds,                                \\
S_2&=&\sum_k\sum_{1\leq i <l}  E \int_{t_i}^{t_{i+1}}
\big[|F_k(s)|_n^2 -| B_{k,i}(\Pi_n
u(t_{i}))  - B_{k,i}( u_{i})|_n^2\big]  \, ds  ,                           \\  
S_3&= & \sum_k E \int_0^{\tau} |F_k(s)|_n^2\, ds.                     
\end{eqnarray*}
For any $\varepsilon>0$ 
\begin{equation*}                                    \label{S1}
S_1\leq \varepsilon\sum_{1\leq i\leq l}  
 \tau \, E\|e_i\|_n^2 
+\tfrac{1}{\varepsilon}R,                                               
\end{equation*} 
where
\begin{align}
R=R(t_l)=&\sum_{1\leq i\leq l}
E\int_{t_{i-1}}^{t_i}
|\Pi_nA(s,u(s))-A_i(\Pi_nu(t_i))|^2_{n^{\ast}}\,ds 
\leq 3  \sum_{1\leq j\leq 3} R_j, 
                                                                 \label{13.08.04}\\
R_1=&\sum_{1\leq i\leq l} E\int_{t_{i-1}}^{t_i}
|\Pi_n A(s,u(s))-\Pi_nA(t_i,u(s))|^2_{n^{\ast}} \,ds   ,                   \nonumber\\  
R_2=&\sum_{1\leq i\leq l} E\int_{t_{i-1}}^{t_i}
|\Pi_n A(t_i,u(s))-A_i(\Pi_nu(s))  |_{n^{\ast}}^2 \, ds ,               \nonumber\\  
R_3=& \sum_{1\leq i\leq l} E\int_{t_{i-1}}^{t_i}
|A_{i}(\Pi_nu(s))-A_i(\Pi_nu(t_i))  |_{n^{\ast}}^2                   
                                                                 \nonumber
\,ds.                                                                         
\end{align}
Due to  condition 
\eqref{HA} on the time regularity of $A$ in 
Assumption \ref{assumption R}, 
  \eqref{bpin}, {\bf (Cn$\tau$)},    
the Lipschitz condition \eqref{LipAAspi}  in Assumption
\ref{assumption ST} and  
 inequality  \eqref{Holderu} from Remark \ref{19.02.04},  
we deduce 
\begin{eqnarray}                                                      \label{majoR1spi}
R_1 &\leq &  
\tau^{2\nu}\,   p^2  
E\int_0^T(K    |u(s)|_{\mathcal V}^2   +\eta)\,ds
,                      \\                  
R_2 & \leq 
& c(\tau^{2\nu}+\varepsilon_n^2)
\Big( E\int_0^T|u(s)|_{\mathcal V}^2\,ds
+ \sum_{1\leq i\leq m} \tau\, E\xi_i^{n,\tau}\Big),                 \\                       
R_3&\leq& L_1p^2\sum_{1\leq i\leq l}\int_{t_{i-1}}^{t_i}
E\|u(s)-u(t_i)\|^2\,ds
\leq 
TL_1p^2M_1\tau , .                                           
\label{14.08.04}
\end{eqnarray}
with $M_1:=C(r_1+r_2+M)$. 
By  \eqref{bpin}, the 
regularity condition {\bf(R3)} on $B$  from Assumption \ref{assumption R},  
the growth condition {\bf(ST2)} on $B_{i,k}$ 
 from Assumption \ref{assumption ST},      
and   
by condition \eqref{Inm}   
on the initial values we have   
\begin{eqnarray}
S_3&\leq&  2 \sum_k\int_0^{\tau}E|\Pi_nB_k(s,u(s))|^2_{n}\,ds
  +2\sum_k{\tau} 
E|B_{k,0}(u_0^{n,\tau})|^2_{n}                  \nonumber \\
&\leq& 2\tau   p^2  
\Big( K\sup_{t\in[0,T]}E\|u(t)\|^2_{\mathcal H}
+\sup_{t\in[0,T]}E\eta(t)\Big)                                  \nonumber\\
                                                          \label{14.06.04}
&& + 2 \tau   \Big(K\sup_{n,m}E\|u^{n,\tau}_0\|_n^2
+\sup_{n,m}Eg_0^{n,\tau}\Big). 
\end{eqnarray}
Using the simple inequality 
$|b|_n^2-|a|_n^2
\leq \varepsilon |a|_n^2
+(1+\tfrac{1}{\varepsilon})|b-a|_n^2$ 
 with 
$$
a:=B_{k,i}(\Pi_nu(t_i))-B_{k,i}(u_i), 
\quad b:=F_k(s), 
$$
for any $\varepsilon>0$ we have 
$S_2\leq \varepsilon P_1+(1+\tfrac{1}{\varepsilon})P_2$ 
with 
\begin{eqnarray}
P_1=P_{1l}&=& \sum_{1\leq i \leq l}  
E\sum_k  \tau | B_{k,i}(\Pi_n u(t_{i})) 
- B_{k,i}( u_{i})|_n^2,       \nonumber \\
                            \label{12.08.04}
P_2=P_{2l}&=&\sum_{1\leq i <l} 
E \int_{t_i}^{t_{i+1}}
\sum_k|\Pi_nB_k(s,u(s))- B_{k,i}(\Pi_n u(t_i))|_n^2 
\, ds. 
\end{eqnarray}
By  Remark \ref{remark Lip} on the Lipschitz
continuity of $B_{k,i}$  we 
get 
$
P_1\leq L_2E \sum_{1\leq i\leq l} \tau \|e_i\|^2_n.
$
Clearly, $P_2\leq 3(Q_1+Q_2+Q_3)$ with 
\begin{eqnarray*}
Q_1 & =&\sum_{1\leq i <l} 
E \int_{t_i}^{t_{i+1}}
\sum_k|\Pi_nB_k(s,u(s))- \Pi_nB_k(t_i,u(s))|_n^2 
\, ds, \\
Q_2& =& \sum_{1\leq i <l} 
E \int_{t_i}^{t_{i+1}}
\sum_k|\Pi_nB_k(t_i,u(s))- B_{k,i}(\Pi_n u(s))|_n^2 
\, ds, \\  
Q_3& =& \sum_{1\leq i <l} 
E \int_{t_i}^{t_{i+1}}
\sum_k
|B_{k,i}(\Pi_n u(s))- B_{k,i}(\Pi_n u(t_i))|_n^2 
\, ds. 
\end{eqnarray*}
Due to 
{\bf (R4) (ii)} in Assumption \ref{assumption R} 
on 
the time regularity of $B$, 
consistency {\bf(Cn$\tau$)}, the  Lipschitz 
continuity of $B_{k,i}$  
proved in Remark \ref{remark Lip}\,, 
\eqref{Holderu}  proved in Remark \ref{19.02.04} 
and \eqref{bpin}, 
\begin{eqnarray*}
Q_1 &\leq &\tau^{2\nu}   p^2 
\Big( K\, E\int_0^T |u(s)|_{\mathcal V}^2   \,ds + T
E\eta\Big),\\  Q_2 &\leq &  c  (\tau^{2\nu}+\varepsilon_n^2)
\Big( E\int_0^T |u(s)|_{\mathcal V}^2  \,ds
+\sup_{n,m}\sum_{0\leq i<l} \tau E\eta_i^{n,\tau}\Big),\\ 
Q_3 &\leq & L_2p^2T\sup_{|t-s|\leq \tau}
E\|u(t)-u(s)\|^2\leq \tau L_2p^2TM_1 .
\end{eqnarray*}
Hence  

\begin{eqnarray}
S_2&\leq&  
\varepsilon L_2  \, E\sum_{1\leq i\leq l} \tau\, 
\|e_i\|^2_n 
+C\big(1+\frac{1}{\varepsilon}\big)(\tau^{2\nu}
+\varepsilon_n^2)                                        \nonumber\\
                                                       \label{12.07.04}
&&
\times \Big(E\int_0^T \!|u(s)|_{\mathcal V}^2  \,ds
 + T E\eta +  \sup_{n,m}\sum_{i}\tau E\eta_i^{n,\tau}
+M_1 \Big), 
\end{eqnarray} 
where $C=C(p,K,L_2,c)$. 
Choosing $\varepsilon>0$ sufficiently small, 
from \eqref{majoE1spi} and \eqref{majoR1spi}--\eqref{12.07.04} 
we obtain for $l=1, \cdots, m$,  
\begin{align}
E|e_l|_n^2 + &
\tfrac{\lambda}{2} E \sum_{1\leq i\leq l}
\tau  \|e_i\|_n^2 \nonumber \\
                                               \label{1estimEspi}
& \leq E|e_0|_n^2   +L\sum_{1\leq i\leq l} 
\tau \, E|e_i|_n^2  
+  C( \tau^{2\nu}+ \varepsilon_n^2)
(r_1+r_2+M), 
\end{align}
where $C=C(K,\lambda, p,T,c,L_1,L_2)$ is a constant. Since
$\sup_m\sum_{i=1}^m\tau=T<+\infty$, if $L\tau < 1$ 
a discrete version of Gronwall's lemma yields
the existence of constants $C_1=C_1(L,\lambda,T)$ 
and 
$C_2=C_2(L,K,\lambda,p,T,c,L_1,L_2)$ such that for sufficiently 
large $m$
$$
\max_{1\leq l\leq m}  E|e_l|_n^2
\leq C_1 E|e_0|_n^2   
+ C_2(r_1+r_2+M)( \tau^{2\nu}+ \varepsilon_n^2)
$$
holds for all $n$. 
This together with \eqref{1estimEspi} concludes the proof of
\eqref{rate1ntau}.
To prove \eqref{ratentau} notice that 
from 
\eqref{diffspi} by the same calculations as above, but 
taking first $\max$ in $l$ and then expectation, 
we get  
$$
E\max_{1\leq i\leq m}|e_i|_n^2
\leq E|e_0|_n^2 
+E\sum_{1\leq i\leq m} \! \tau \|e_i\|_n^2 + 
L\sum_{1\leq i\leq l} \tau  E|e_{i}|_n^2
$$
\begin{equation}                             \label{majodiffspi}  
+R(T)
+ R_0 
+E\max_{0\leq i\leq m}I(t_i)+EQ(T).  
\end{equation}
where  $R(T)=R(t_m)$ is defined by \eqref{13.08.04} and    
\begin{equation}                                \label{noR0}
R_0=P_{1m}+2P_{2m}+S_3. 
\end{equation}
The terms $R(T)$, $P_{1m}$ and $P_{2m}$ have already been 
estimated above by the right-hand side of  
\eqref{1estimEspi} and $S_3$ has been estimated by
\eqref{14.06.04}. Notice that   
$$ 
EQ(T) =E  \int_0^T\sum_k |F_k(s)|_{n}^2\,ds\leq 2P_{1m}+2P_{2m}, 
$$
and by Davis' inequality 
\begin{eqnarray*}
E\max_{0\leq i\leq m} I(t_i)&\leq &6E\left\{\int_{0}^{T} 
\sum_k|(e(s)\, ,\,F_k(s))_n|^2\,ds\right\}^{1/2}         \nonumber \\
                                                  \label{15.11.04}
& \leq & \tfrac{1}{2}E\max_{0\leq i<m}|e_i|^2_n
+18E\int_0^T \sum_k |F_k(s)|^2_n\,ds.    
\end{eqnarray*} 
Thus from  
\eqref{majodiffspi}
we obtain \eqref{ratentau}. 
\end{proof}
\smallskip

\begin{rk}
One can show, like it is observed in \cite{GyMi2},  
that  if instead of the Lipschitz condition 
\eqref{LipAAspi} we assume that $A^{n,\tau}_i$ 
are hemicontinuous and $B^{n,\tau}_{k,i}$  
satisfy the Lipschitz condition \eqref{lipBnm},    
then the order of the 
speed of convergence is divided by two.
\end{rk}

\subsection{Examples}                         \label{exspi}
(i)  {\rm Consider from Example \ref{ExSpace} 
the normal triples 
$$
V\hookrightarrow H^*\hookrightarrow V^*, 
\quad V_n\hookrightarrow H_n^*\hookrightarrow V_n^*
$$
with the orthogonal projection 
$\Pi_n: H=L^2(\mathbb R^d)\to H_n$ 
and auxiliary spaces  
$\mathcal H=W^{r+\rho}_2(\mathbb R^d)$ and 
$\mathcal V=W^{r+l}_2(\mathbb R^d)$ for some 
$l>\rho\geq0$. 

Let  $A$ and   $B=(B_k)$ be 
$\mathcal P\otimes\mathcal B(V)$-measurable mappings 
from $[0,\infty[\times \Omega\times V$ 
into $V^*$ and $H^{d_1}$, respectively, 
satisfying Assumptions \ref{E} and 
\ref{assumption R} such that $f$ and $g$ in 
\eqref{R02} satisfies 
$$
\sup_{t\in[0,T]}f(t)\leq M, \quad \sup_{t\in[0,T]}g(t)\leq M. 
$$ 
 For $\omega\in \Omega$, $j=1,\dots,m$ and $i=0,\dots,m-1$  
let $A^{n,\tau}_j(\omega,\cdot):V^n\to V_n^*$ and 
$B^{n,\tau}_{k,i}(\omega,\cdot):V^n\to H_n$ 
be defined 
by 
\begin{equation}                             \label{choice1}
\<A^{n,\tau}_j(\omega,u)\,,v\>_n
=\<A(t_j,\omega,u,)\, ,v\> 
\quad \text {and $B^{n,\tau}_{k,i}(\omega,u)=\Pi_n
B_k(t_i,\omega,u)$}
\end{equation}
for all $u,v\in V_n$, where $\<\,,\>_n$ denotes 
the duality between $V_n$ and $V_n^*$.  
Then it is easy to see, like in 
Example \ref{ExSpace}, that due to 
{\bf (1)}, {\bf (2)} and {\bf (3)} 
in Assumption \ref{E}, 
{\bf (ST1)}, {\bf (ST2)} and {\bf (ST3)} 
in Assumption \ref{assumption ST} hold respectively.  
In the same way 
as {\bf (Cn)} is verified in Example \ref{ExSpace}, 
one can also
easily show that the consistency  assumption {\bf (Cn$\tau$)}
holds.  

(ii) Another  choice for $A^{n\tau}_j$ 
and $B^{n\tau}_{k,i}$ can be defined by 
$$
\<A^{n,\tau}_{j}(u)\,,v\>_n=\frac{1}{\tau}
  \int_{t_{j-1}}^{t_j} \<A(s,u)\,,v\> ds, 
\quad u,v\in V_n,
$$
\begin{equation}                           \label{choice2}
\quad B^{n,\tau}_{k,0}(u)=\Pi_n B_{k}(0,u),
\quad 
B^{n,\tau}_{k,j}(u)=\frac{1}{\tau}  
\int_{t_{j-1}}^{t_j} \Pi_n B_k(s,u)\, ds, 
\quad u\in V_n, 
\end{equation}
instead of \eqref{choice1}. 
One can show by a similar computation as before, 
combined with the use of  Jensen's inequality, that  
Assumption \ref{E} and {\bf(R3)}-{\bf(R4)} 
in Assumption \ref{assumption R} imply  
Assumption \ref{assumption ST}
and  condition {\bf (Cn$\tau$)}.

(iii) Finally,  let $V_n=V$, $H_n=H$ and 
let $\Pi_n$ be the identity operator for every $n$. 
Let 
Assumptions \ref{E} and 
\ref{assumption R} 
hold. Then 
one recovers the conclusions 
of Theorems 3.2 and 3.4 in
\cite{GyMi2} concerning the
rate of convergence of the implicit time discretization scheme
with $\varepsilon_n=0$. 

\smallskip

\section
{Explicit space-time discretization scheme}              \label{esp}
\subsection{Description  of the scheme}

Let $V_n$, $H_n$ and $V_n^*$ be a normal
triple and $\Pi_n$ be continuous
linear operators which satisfy the
condition \eqref{bpin}. 
Assume moreover that for each 
$n \geq 0$ 
as sets 
$$ V_n=H_n=V_n^{\ast},  $$
and there is a constant  $\vartheta(n)$ 
such that
\begin{equation}                                           \label{VnHn}
\|u\|_{V_n}^2 \leq \vartheta(n)\, |u|^2_{H_n}\; , 
\quad \forall u\in H_n.
\end{equation}
Then by duality we also have 
\begin{equation*}                                         \label{HnVn*}
|u|_{H_n}^2 \leq \vartheta(n)\, |u|^2_{V_n^*}\; ,
\quad \forall u\in V_n^*.
\end{equation*}
Consider for each $n$ and $i=0,1,\cdots,m-1$ 
the equations 
\begin{equation}                                         \label{defexp}
u^{n}_{\tau,i+1}= u^{n}_{\tau,i}
+ \tau\,  A^{n,\tau}_{i}
\big(u^{n}_{\tau,i}\big)                            
+\sum_k {B}^{n,\tau}_{k,i} 
\big(u^{n}_{\tau,i}\big)\,
\big(W^k(t_{i+1})- W^k({t_i})\big),  
\end{equation}
for $V_n$-valued $\mathcal F_{t_i}$-measurable 
random variables 
$u^{n}_{\tau,i}$ for $i=1,\cdots,m$,  
where $u_{\tau,0}^{n}$ is a given $V_n$-valued 
 $\mathcal F_0$-measurable random variable, and  
$$
A^{n,\tau}_{i} : \Omega \times V_n  \to V_n^*\;   
\quad \mbox{\rm and} 
\quad  B_{k,i}^{n,\tau} : \Omega\times V_n 
\to H_n 
$$
are given 
${\mathcal F}_{t_i}\otimes {\mathcal B}(V_n)$-measurable  
mappings such that Assumption \ref{assumption ST}
holds.

\begin{proposition}                           \label{10.12.04}
Let Assumption \ref{assumption ST}
hold. 
Then for any $V$-valued $\mathcal F_0$-measurable 
random variable $u^n_{\tau,0}$ such that 
$E\|u^n_{\tau, 0}\|^2_{V_n}<\infty$,  
the system of equations 
\eqref{defexp} has 
a unique solution $(u^n_{\tau,i})_{i=1}^m$ 
such that $u^n_{\tau,i}$ is 
$\mathcal F_{t_i}$-measurable and 
$E\|u^n_{\tau,i}\|^2_{V_n}<\infty$ for all $i$, $m$ and $n$. 

\begin{proof} By \eqref{VnHn} we have 
$\|u^n_{\tau,i+1}\|_n^2
\leq \vartheta(n)|u^n_{\tau,i+1}|^2_n$, 
and by \eqref{defexp} 
$$
E|u^n_{\tau,i+1}|_{n}^2\leq 3E|u^n_{\tau,i}|^2_n
+3\tau E|A^{n,\tau}_i(u^n_{\tau,i})|^2_n
+3\tau\sum_k E|B^{n\tau}_i(u^n_{\tau,i})|^2_n
$$
$$
\leq 3
\big(\vartheta(n)+\vartheta(n)\tau K+\tau K\big)\, 
E\|u^n_{\tau,i}\|^2_n
+3\tau\big(\vartheta(n)+1\big)M.
$$
Hence we get the proposition by induction on $i$. 
\end{proof}
\end{proposition}

\subsection{Rate of convergence of the  scheme}

The following theorem gives the 
rate of convergence 
of $e^n_{\tau,i}:=\Pi_n u(t_i)-u^n_{\tau,i}$.

\begin{Th}                                                 \label{thexp}
Let  
Assumption \ref{assumption R}, Assumption 
\ref{assumption ST} with index $j=i=0,\dots, m-1$ in its
formulation, and the consistency  condition 
{\bf (Cn$\tau$)} hold. Let $n$ and $\tau$ satisfy 
\begin{equation}                                        \label{stability}
L_1\tau\vartheta(n) 
+2\sqrt{L_1L_2\tau\vartheta(n)}\leq q
\end{equation} 
for some constant $q<\lambda$, 
where $L_1$ 
and $L_2$ are the Lipschitz constants 
in \eqref{LipAAspi} and \eqref{lipBnm}, 
respectively.  Then 
\begin{equation}                                            \label{rateexp}
 E\max_{0\leq i\leq m}   
|e^{n}_{\tau,i}|_{H_n}^2 +
\sum_{0\leq i<m} \tau E\|e^n_{\tau,i}\|_{V_n}^2           
\leq C_1  E |e^n_{\tau,0}|_{H_n}^2    
 + C_2( \tau^{2\nu}+ \varepsilon_n^2)
(r_1+r_2+M), 
\end{equation} 
where $C_1=C_1(\lambda, q,L,T)$ and 
$C_2=C_2(\lambda,q, L,K,T,p,c,L_1,L_2)$ are constants.
\end{Th}

\begin{proof} Note that when we refer 
to any condition in Assumption \ref{assumption ST} 
then we mean it with the index $j$ replaced 
in its formulation with 
$i$ running through $0,\cdots,m-1$. To ease notation
we omit  the indices $n$ and $\tau$ from 
$e^n_{\tau,i}$, $u^n_{\tau,i}$, $A^{n,\tau}_i$ and $B^{n,\tau}_i$ 
when this does not cause ambiguity. 
For any $i=0,\, \cdots,\, m-1$ 
\begin{align*}
|e_{i+1}|_n^2 & 
- |e_i|_n^2 
=  2\, \int_{t_i}^{t_{i+1}}
\langle e_i\, ,\,  \Pi_n  A(s,  u(s))
- A_i(u_i))\rangle_n \, ds \\
& + 2\,\sum_k \int_{t_i}^{t_{i+1}} 
(e_i\, ,\,  F_k(s))_n \,dW^k(s)
+  \sum_k  \left| 
\int_{t_i}^{t_{i+1}}F_k(s)\, dW^k(s) 
\right|_n^2\\
& +  \Big|\int_{t_i}^{t_{i+1}} 
\big[ \Pi_n A(s,  u(s)) -
A_i(u_i) \big] \, ds \Big|^2_n \\
& +2\sum_k  \Big(
 \int_{t_i}^{t_{i+1}} 
\big[ \Pi_n A(s,u(s)) -A_i(u_i)\big] \, ds \, ,
 \int_{t_i}^{t_{i+1}} F_k(s)\,dW^k(s) 
\Big)_n , 
\end{align*}
where 
$$                                          
 F_k(s)
=\Pi_n  B_k(s,u(s)) - {B}_{k,i}(u_i), 
\quad s\in]t_i,t_{i+1}], \quad i=0,1,\cdots,m-1. 
$$
Hence for $l=1, \cdots, m$ and every $\delta >0$,  
\begin{eqnarray}
|e_l|^2_n&\leq & |e_0|^2_n
+2\sum_{0\leq i<l}
\int_{t_i}^{t_{i+1}}
\langle e_i\, ,\,  \Pi_n  A(s,  u(s))
- A_i(u_i))\rangle_n \, ds
+2\, I(t_l)+Q(t_l)                                 \nonumber \\
&&+(1+\tfrac{1}{\delta})S(t_l)+\delta Q(t_l) ,     \label{eq5.6}
\end{eqnarray}
 where 
\begin{eqnarray*}
I(t_l)&=&\int_0^{t_l}(e(s) \,,F_k(s))\,dW^k(s), 
\quad
e(s) :=e_i\quad\text{for $s\in]t_i,t_{i+1}]$, $i\geq 0$}, \\
S(t_l)&=&\sum_{0\leq i<l}
\Big|\int_{t_i}^{t_{i+1}} 
\big[ \Pi_n  A(s,  u(s)) -A_i(u_i)\big]\, ds \Big|^2_n, \\
Q(t_l)&=&\sum_{0\leq i<l}
\sum_k 
\left| \int_{t_i}^{t_{i+1}}F_k(s)\, dW^k(s) \right|_n^2. 
\end{eqnarray*}
First we prove 
\begin{equation}                                  \label{rate1exp}  
\max_{1\leq i\leq m}
 E |e_i|_{n}^2 +
\sum_{0\leq i<m} \tau  E\|e_i\|_{n}^2         
\leq C_1 E|e_0|_{n}^2    
+C_2( \tau^{2\nu}+ \varepsilon_n^2)
(r_1+r_2+M), 
\end{equation} 
with some constants 
$C_1=C_1(\lambda, q,L,T)$ and 
$C_2=C_2(\lambda,q, L,K,T,p,c,L_1,L_2)$. 
To this end we take expectation in both sides of 
\eqref{eq5.6} and use  
the strong monotonicity condition 
\eqref{monotonspi} in Assumption \ref{assumption ST}, 
to get  
\begin{align*}
E|e_l|^2_n &  \leq  E|e_0|^2_n  +2 E\sum_{0\leq i<l}
\tau \langle e_i\, ,
 \,A_i(\Pi_n u(t_i))-A_i(u_i)\rangle_n              \\  
&\quad  +E\sum_{0\leq i<l}
\tau |B_{k,i} (\Pi_n u(t_i))-B_{k,i} (u_i)|^2_n
+\sum_{i=1,2} S_i 
+(1+\tfrac{1}{\delta})ES(t_l)+\delta EQ(t_l)     \\ 
\leq&  E|e_0|^2_n  
-\lambda E\sum_{0\leq i<l} \! \tau \|e_i\|^2_n
+L E\sum_{0\leq i<l}\!  \tau |e_i|^2_n          
 +\sum_{i=1,2} \! S_i
+(1+\tfrac{1}{\delta})ES(t_l)+\delta EQ(t_l), 
\end{align*}
for any $\delta>0$, where 
\begin{eqnarray*}
S_1&=&2\,\sum_{0\leq i<  l} 
E\int_{t_{i}}^{t_{i+1}}
\langle e_i\, ,\, \Pi_n 
A(s,u(s))-A_i(\Pi_n u(t_i))\rangle_n\,ds,                                \\
S_2&=&\sum_k\sum_{0\leq i <l}  E \int_{t_i}^{t_{i+1}}
\big[ |F_k(s)|_n^2 -| B_{k,i}(\Pi_n
u(t_{i}))  - B_{k,i}( u_{i})|_n^2\big] \, ds.                                  
\end{eqnarray*}
As in the proof of Theorem \ref{th1spi} 
we get  for any $\varepsilon >0$, 
\begin{eqnarray*}
S_1 &\leq & \varepsilon  \sum_{0\leq i < l}  
\tau E\|e_i\|_n^2
+\tfrac{1}{\varepsilon}
C(r_1+r_2  +M  )(\tau^{2\nu}+\varepsilon_n^2 ) ,\\
S_2&\leq &L_2 \, \varepsilon \sum_{0\leq i < l} 
\tau  E\|e_i\|_n^2 + \tfrac{1}{\varepsilon}
C(r_1+r_2 +M  )(\tau^{2\nu}+ \varepsilon_n^2  ) 
\end{eqnarray*}
with a constant  $C=C(K,p,T,L_1,L_2,c)$.  
Notice that for any $\varepsilon >0$,  
\begin{align}
ES(t_l)\leq & \tau \vartheta(n) J(t_l),            \nonumber \\
                                                   \label{15.10.04}
J(t_l):=& \sum_{0\leq i<l}
E\int_{t_i}^{t_{i+1}}
|\Pi_n A(s,u(s))-A_i(u_i)|_{n^{\ast}}^2\,ds
\leq (1+\varepsilon)R_0
+(1+\tfrac{1}{\varepsilon})R, \\
                                                 \label{16.10.04}
EQ(t_l)\leq & (1+\varepsilon)P_1
+(1+\tfrac{1}{\varepsilon})  P_2 , 
\end{align}
 where        
\begin{eqnarray}                                 \label{19.10.4} 
P_1&:=& \sum_{0\leq i <l}
E\sum_k| B_{k,i}
(\Pi_n u(t_{i})) - B_{k,i}( u_{i})|_n^2\tau
\leq L_2E\sum_{0\leq i <l} \tau \|e_i\|_n^2      \nonumber\\
P_2& :=& \sum_{0\leq i <l} 
E \int_{t_i}^{t_{i+1}}
\sum_k|\Pi_n B_k(s,u(s))- B_{k,i}(\Pi_n u(t_i))|_n^2 
\, ds,                                            \nonumber\\
                                                  \nonumber
R_0& :=& E\sum_{0\leq i<l}   \tau 
|A_i(\Pi_n u(t_i))-A_i(u_i)|^2_{n^{\ast}}
\leq L_1 E\sum_{0\leq i<l}
\tau \|e_i\|^2_{n},                               \\
R&:=&E\sum_{0\leq i<l}\int_{t_i}^{t_{i+1}}
|\Pi_n A(s,u(s))-A_i(\Pi_n u(t_i))|^2_{n^{\ast}}\,ds,  
                                                  \nonumber
\end{eqnarray}
for any $l=1,2,\cdots,m$.
In the same way as in the proof 
of Theorem \ref{th1spi} 
we obtain 
\begin{equation}                                  \label{18.10.04}
R\leq C(\tau^{2\nu}+\varepsilon_n^2)(r_1+r_2+M), 
\end{equation}
and 
that 
\begin{equation}                                  \label{20.10.04}
P_2\leq 
C^{\prime}(\tau^{2\nu}+\varepsilon_n^2)
(r_1+r_2+M), 
\end{equation}
where $C=C(K,p,c,L_1,T)$  and 
$C^{\prime}=C^{\prime}(K,p,c,L_2,T)$ are constants. 
Consequently, 
\begin{eqnarray}
E|e_l|^2_n& \leq  &
E|e_0|^2_n   
+ (\mu-\lambda) E\sum_{0\leq i<l} \tau \|e_i\|^2_n
+LE\sum_{0\leq i<l} \tau |e_i|^2_n                 \nonumber \\ 
                                                \label{10.010.07}
&& +\big(1+\tau \vartheta(n)\big)
(1+\tfrac{1}{\delta}+\tfrac{1}{\varepsilon})
C(\tau^{2\nu}+\varepsilon_n^2)(r_1+r_2+M) 
\end{eqnarray}
for any $\delta>0$ and $\varepsilon>0$, 
where 
$$
\mu=(1+\varepsilon)\big[ (1+\tfrac{1}{\delta})
\tau \vartheta(n) L_1+\delta L_2\big]  +\varepsilon(1+L_2),  
$$
and $C=C(K,p,c,T,L_1,L_2)$ is a constant. It is easy to see 
that due to \eqref{stability}
$$
\inf_{\delta>0}
(1+\tfrac{1}{\delta})
\tau \vartheta(n) L_1+\delta L_2  
=\tau\vartheta(n)L_1 +2\sqrt{\tau\vartheta(n)L_1L_2}\leq q. 
$$
Therefore we can take $\delta>0$ and $\varepsilon>0$ 
such that 
$
\mu\leq (q+\lambda)/2.    
$
Thus 
from \eqref{10.010.07} we can get 
\begin{eqnarray*}
E|e_l|^2_n &
\leq &
E|e_0|^2_n 
-\tfrac{1}{2}(\lambda-q) E\sum_{0\leq i<l}  \tau \|e_i\|^2_n
+L E\sum_{0\leq i<l} \tau |e_i|^2_n \\
&& + C(\tau^{2\nu}+\varepsilon_n^2)(r_1+r_2+M),  
\end{eqnarray*}
with a constant $C=C(K, \lambda,q,p,c,T,L_1,L_2)$.   
Hence by a discrete version of Gronwall's lemma 
we obtain \eqref{rate1exp}.  
To prove \eqref{rateexp} note that 
\eqref{eq5.6}  yields 
\begin{equation}                           \label{21.04.07}
E\max_{1\leq l\leq m}|e_l|^2_n
\leq  |e_0|^2_n  +E\sum_{0\leq i<l}
\tau  \|e_i\|^2_n+ 2 E S(T) 
+ 2  E\max_{1\leq l\leq m}  I(t_l)+2EQ(T),  
\end{equation}
where by \eqref{15.10.04}--\eqref{18.10.04}  
 $E S(T)\leq  \tau \vartheta(n) J(T)$, 
 and   
\begin{eqnarray*}
J(T)&\leq& 2L_1E\sum_{0\leq i<m}  \tau \|e_i\|^2_n
+2C(\tau^{2\nu}+\varepsilon_n^2)(r_1+r_2+M). 
\end{eqnarray*}
By \eqref{16.10.04}, \eqref{19.10.4} and 
\eqref{20.10.04}
$$
EQ(T)\leq 2L_2\sum_{0\leq i<m}  \tau  \|e_i\|^2_n
+2C^{\prime}(\tau^{2\nu}+\varepsilon_n^2)(r_1+r_2+M). 
$$ 
Finally,  in the same way as 
 equation \eqref{15.11.04} is obtained, we get  
\begin{align*}
E \max_{1\leq i\leq m} I(t_i) &\leq 6E\left\{\int_{0}^{T} 
\sum_k|(e(s)\, ,\,F_k(s))_n|^2\,ds\right\}^{1/2}   \\
& \leq \tfrac{1}{2}E\max_{0\leq i<m}|e_i|^2_n
+18 E  \int_0^TF_k^2(s)\,ds\leq 
\tfrac{1}{2}E\max_{0\leq i<m}|e_i|^2_n
+18EQ(T).   
\end{align*} 
Consequently, from \eqref{21.04.07} we 
obtain \eqref{rateexp} by \eqref{rate1exp}. 
\end{proof}

\subsection{Example}                              \label{example Enm} 
Consider
again the spaces 
$$
\mathcal V\subset\mathcal H
\subset V\hookrightarrow H\hookrightarrow V^*, 
\quad V_n\hookrightarrow H_n \hookrightarrow V_n^*
$$
from Examples \ref{ExSpace} and \ref{exspi}. 
Notice that $V_n$, $H_n$ and $V_n^*$ are 
identified as sets 
and that due to the converse 
inequality (\ref{converse}) we have 
\eqref {VnHn} with $\vartheta(n)=C2^{nr}$. 
  
Let $A$ and $B$ satisfy the 
same conditions as in Example \ref{exspi}. 
Define $A^{n,\tau}_j$ and $B^{n,\tau}_{k,i}$ 
for $j=i=0,\dots,m-1$  
by \eqref{choice1} or define 
$A_j^{n,\tau}$ 
by \eqref{choice1} for $j=0$ and 
$A_j^{n,\tau}$, $B^{n,\tau}_{k,i}$ 
by \eqref{choice2} for $j=1,\dots,m-1$ 
and 
$i=0,\dots,m-1$.    
Then, as shown in section \ref{exspi},
$A^{n,\tau}_j$ and $B^{n,\tau}_{k,i}$ 
 satisfy the conditions 
in Assumption \ref{assumption ST}  
as well as {\bf (Cn$\tau$)}.   
Hence, if the solution $u$ satisfies 
{\bf (R1)}-{\bf(R2)} in Assumption 
\ref{assumption R},  
and  
$L_1 \tau \vartheta_n 
+2\sqrt{L_1 L_2 \tau \vartheta_n}\leq q<\lambda$, 
then 
the conditions of Theorem \ref{thexp} hold.

\section{Examples of approximations of  stochastic PDEs}

In this section we present some examples 
of stochastic PDEs 
for which the previous theorems provide 
rates of convergence for the
above space and space-time discretization schemes. 
We refer to section 5 in  \cite{GyMi2} for more details.
In this section for integers $l$ the notation 
$|u|_l=|u|_{W^l}$ 
means the norm of $u$ in $H^l=W^l(\mathbb R^d)$. 

\subsection{Quasilinear equations}
\noindent Let us consider the stochastic
partial differential equation
\begin{align}
du(t,x)=& \big( Lu(t)
+F(t,x,\nabla u(t,x),u(t,x)\big)\,dt
                                                   \nonumber
\\
                                                   \label{ex1eqnl}
&+\sum_k
\big( M_ku(t,x)+g_k(t,x)\big)\,dW^k(t),
\quad t\in(0,T], \,x\in{\mathbb R}^d, 
\end{align}
with initial condition
\begin{equation}                                   \label{ini}
 u(0,x)=u_0(x),
\quad x\in\mathbb R^d, 
\end{equation}
where $F$ and $g_k$ are Borel
functions of
$(\omega,t,x,p,r)\in \Omega\times[0,\infty)
\times{\mathbb R}^d\times{\mathbb R}^d
\times{\mathbb R}$  and of $(\omega,t,x)\in
\Omega\times[0,\infty)\times{\mathbb R}^d$,
respectively, and
$L$,
$M_k$ are differential operators of the  form
\begin{equation*}                                      \label{linop}
L(t)v(x)=\sum_{|\alpha|\leq 1, |\beta|\leq1} \!\!
D^{\alpha}(a^{\alpha\beta}(t,x)D^{\beta}v(x)),
\;  M_k(t)v(x)=\sum_{|\alpha|\leq1}
b^{\alpha}_k(t,x)D^{\alpha}v(x),
\end{equation*}
with functions
$a^{\alpha\beta}$ and $b^{\alpha}_k$ of
$(\omega,t,x)\in \Omega\times[0,\infty)\times{\mathbb R}^d$,
for all multi-indices
$\alpha=(\alpha_1,...,\alpha_d)$,
$\beta=(\beta_1,...,\beta_d)$ of length
$|\alpha|=\sum_i\alpha_i\leq 1$, $|\beta|\leq1$.
Here, and later  on
$D^{\alpha}$ denotes $D^{\alpha_1}_1...D^{\alpha_d}_d$
for any multi-indices
$\alpha=(\alpha_1,...,\alpha_d)\in\{0,1,2,...\}^d$,
where
$
D_i=\frac{\partial}{\partial x_i}
$
and $D_i^0$ is the identity operator.
We use the notation
$\nabla_p:=(\partial/\partial p_1,...,\partial/\partial p_d)$.

Let $K$ and $M$ denote some 
non-negative numbers. Fix an integer $l\geq 0$ and suppose 
that the following conditions hold:

\noindent
{\bf Assumption (A1)}
(Stochastic parabolicity).
{\it
There exists a constant
$\lambda>0$ such that
\begin{equation*}                                   \label{22.29.05.06}
\sum_{|\alpha| = 1,|\beta| = 1} 
\left(
a^{\alpha\beta}(t,x)-\tfrac{1}{2} 
\sum_k \big( b^{\alpha}_k
b^{\beta}_k\big) (t,x)
\right) 
\,z^{\alpha}\, z^{\beta}\geq \lambda
\sum_{|\alpha|=1}|z^{\alpha}|^2
\end{equation*}
for all $\omega\in\Omega$, $t\in[0,T]$,
$x\in{\mathbb R}^d$ and
$z=(z^1,...,z^d)\in{\mathbb R}^d$, where
$z^{\alpha}
:=z_1^{\alpha_1}z_2^{\alpha_2}...z_d^{\alpha_d}$
for $z\in{\mathbb R}^d$ and multi-indices
$\alpha=(\alpha_1,\alpha_2,...,\alpha_d)$.
}
\smallskip

\noindent
{\bf Assumption (A2)}
(Smoothness of the initial condition).
{\it Let $u_0$ be 
$W^2_l$-valued 
${\mathcal F}_0$-measurable random variable
such that                                   
$  E|u_0|^2_l\leq M$. 
}
\smallskip

\noindent
{\bf Assumption (A3)} 
(Smoothness of the linear term).
{\it  The derivatives of $a^{\alpha\beta}$
and $b^{\alpha}_k$ up to
order $l$  are
${\mathcal P}\otimes{\mathcal B}
(\mathcal{\mathbb R}^d)$
-measurable real functions
such that almost surely 
\begin{equation*}                                    \label{23.29.05.06}
|D^{\gamma}a^{\alpha\beta}(t,x)|
+ |D^{\gamma} b^{\alpha}_k(t,x)|\leq K,
\quad
\text{\rm for all
$|\alpha|\leq1$, $|\beta|\leq1$,
$k=1,\cdots, d_1$,}
\end{equation*}
$t\in[0,T]$,
$x\in{\mathbb R}^d$ and multi-indices
$\gamma$ with $|\gamma|\leq 2$.
}
\smallskip

\noindent
{\bf Assumption (A4)} (Smoothness of the nonlinear term).
{\it  The function $F$ and their
first order partial derivatives in $p$, $x$
and $r$ are
${\mathcal P}\otimes{\mathcal B}(\mathcal{\mathbb R}^d)
\otimes{\mathcal B}(\mathcal{\mathbb R}^d)
\otimes{\mathcal B}(\mathcal{\mathbb R})$-measurable
functions. The function
$g_k$  and  its
derivatives in
$x$  are
${\mathcal P}
\otimes{\mathcal B}(\mathcal{\mathbb R}^d) $-measurable
 functions for every $k=1,..,d_1$.
There exists a constant $K$ and a 
${\mathcal P}\otimes{\mathcal B} $-measurable
 function $\xi$ of $(\omega,t,x)$
such that almost surely 
\begin{align}
& |\nabla_p F(t,x,p,r)| 
+|\tfrac{\partial}{\partial r} F(t,x,p,r)|
\leq {K} \, ,
                                                 \nonumber\\
& |F(t,\cdot,0,0)|_0^2 
+\sum_k  |g_k(t,\cdot)|_2^2  \leq\eta \, ,
                                                 \nonumber\\
& |\nabla_xF(t,x,p,r)|\leq L(|p|+|r|)+\xi(t,x),
\quad |\xi(t)|_0^2\leq\eta 
                                                 \nonumber
\end{align}
for all $t,x, p,r$, where $\eta$ is 
a random variable such that
 $E\eta\leq  M$}.

\noindent

\smallskip

Set $H=L^2({\mathbb R}^d)= W^0_2$,  
$V=W^1_2$,  
${\mathcal H}=W^2_2$ and  ${\mathcal V}=
W^3_2$  and 
 suppose that the assumptions 
{\bf (A1)}--{\bf (A4)} hold with $l=2$. 
Then the operators
$$
A(t,\varphi)= L(t)\varphi
+F(t,.,\nabla \varphi,\varphi),  
\quad 
B_k(t,\varphi)=M_k(t)\varphi + g_k(t,.), 
\quad \varphi\in V
$$
and $u_0$ satisfy the conditions of 
Theorem  \ref{existu}. Hence 
\eqref{ex1eqnl}--\eqref{ini} has 
a unique solution $u$ on $[0,T]$. 
Furthermore, $u$   has a $W^2_2$-valued 
continuous modification 
such that
\begin{equation*}                                       \label{boundW23}
E\sup_{0\leq t\leq T}|u(t)|_2^2
+E\int_0^T|u(t)|^2_3\,dt<\infty . 
\end{equation*}
Consequently  the regularity conditions     
{\bf (R1)} and
{\bf (R2)} in Assumption \ref{assumption R} hold.}  
It is easy to check 
that $A$ and $B_k$ verify condition {\bf (R3)}. 
\smallskip

\noindent
{\bf Assumption (A5)} (Time regularity of $A$ and $B$)
{\it 
Almost surely 

(i)
\begin{align*}
\sum_k |D^\gamma (b^{\alpha}_k(t,x)-b^{\alpha}_k(s,x)) |^2
&\leq K |t-s|,
                                                    \nonumber\\
\sum_k |g_k(s,.)-g_k(t,.)|_1^2 &\leq \eta\, |t-s|.
                                                    \nonumber
\end{align*}
(ii)
 \begin{align*}
|D^\gamma (a^{\alpha,\beta}(t,x)-a^{\alpha,\beta}(s,x)) |^2
&\leq K |t-s|,
                                                   \nonumber\\
|F(t,x,p,r)-F(s, x,p,r)|^2 &  \leq K\, |t-s|\, (|p|^2+|r|^2) , \\
 |\nabla_x F(t,x,p,r)- \nabla_x F(s,x,p,r)|^2 &
\leq K\, |t-s|\, (|p|^2+|r|^2),
                                                    \nonumber\\
 |\nabla_p F(t,x,p,r)- \nabla_p F(s,x,p,r)|^2 
 & \leq K\, |t-s|, \nonumber \\
|\tfrac{\partial}{\partial_r} F(t,x,p,r)- \tfrac{\partial}{\partial_r}
F(s,x,p,r)|^2 & \leq K\, |t-s| . \nonumber
 \end{align*}
for all $|\alpha|\leq1$, $|\beta|\leq1$, $|\gamma|\leq 1$, 
$s,t\in [0,T]$ 
and $x\in {\mathbb R}^d$, 
where $K$ is a constant and  $\eta$ is a random variable 
such that $E\eta\leq M$.}

Clearly, Assumptions (i) and (ii) of {\bf (A5)}   
imply conditions (i) and (ii) 
of 
{\bf (R4)} in Assumption \ref{assumption R}, 
respectively with $\nu=1/2$. 

Let $H_n$, $V_n$ and $\Pi_n$ be defined as in 
Example \ref{ExWavelets} 
and let
$A^n(t,u)$ and $B^n_k(t,u)$ 
be defined by (\ref{exoperatorspace}).
Let $u_0\in W^2_2={\mathcal H}$ and $u_{0}^n =\Pi_n u_0$. 
Recall Example \ref{ExSpace} and notice that 
we can apply Theorem \ref{raten}, and 
by making use of \eqref{direct} we get 
the estimate 
$$
E\sup_{0\leq t\leq T} |u^n(t) -u(t)|_0^2 
+ E\int_0^T |u^n(t)  -u(t)|_1^2\, dt  \leq C\, 2^{-2n},
$$
with a constant $C$ independent of $n$. 
Assume now also {\bf (A5)}, recall 
Example \ref{exspi} and  
define 
$A^{n,\tau}$ and $B^{n,\tau}$ by 
\eqref{choice1}. Notice that 
we can apply Theorem \ref{th1spi}.     
Hence  if $u^{n,\tau}_0=\Pi_n u(0)$   
we get the estimate 
\begin{equation}                             \label{exratespi2}
 E\max_{0\leq i\leq m} 
|u^{n,\tau}_i - u(i\tau)|_0^2
 +\tau E\sum_{0\leq i \leq m} 
|u^{n,\tau}_i - u(i\tau)|_1^2  \leq C\,
\Big(\tau+2^{-2n} \Big).
 \end{equation}
Finally recall Example \ref{example Enm} and 
define $A^{n,\tau}$ and $B^{n,\tau}$ 
as in Example \ref{example Enm}. 
Then we can apply Theorem 
\ref{raten}, and if  
$u^n_{\tau,0}:=\Pi_n u(0)$ and 
${T}2^{2n}/m \leq \gamma$ 
for some constant $\gamma <c\lambda $, 
then we get estimate \eqref{exratespi2} 
for 
the explicit space-time approximations $u^{n}_{\tau,i}$,   
in place of $u^{n,\tau}_i$, with some constant $C$.  
\smallskip

Let us now recall 
Example \ref{ExFiniteEl} and approximate 
\eqref{ex1eqnl}--\eqref{ini} by 
finite difference schemes. Consider 
first the following system of SDEs, 
corresponding to the space discretization 
with finite differences for fixed $h\in(0,1)$:  
\begin{eqnarray}
dv(t)&=&
\big( L_h(t)v(t)
+F_h(t,\nabla_h v(t),v(t))\big)\,dt     \nonumber    \\
                                                      \label{12.28.03}
&&+\sum_{k}\big(M_{k,h}(t)v(t)
+g_{k,h}(t)\big)\,dW^k(t), 
\quad z\in {\mathbb G}=hZ^d ,                               \\
                                                     \label{13.28.03}
v(0)&=&  (u_0(z))_{z\in {\mathbb G}},
\end{eqnarray}
where  $g_{k,h}(t)=(g_k(t,z))_{z\in\mathbb G}$, 
$F_h(t,p,r)=(F(t,z,p,r)_{z\in {\mathbb G}})$ and  
\begin{align}                                       \label{14.28.03}
L_h(t)\varphi :=&\sum_{|\alpha|\leq1,|\beta|\leq1}
\delta_{-}^{\alpha}(a^{\alpha\beta}(t,\cdot)
\delta^{\beta}_+\varphi), 
\quad 
\nabla_h\varphi:=(\delta_{1}\varphi,
\delta_{2}\varphi,\dots,\delta_{d}\varphi),   \\
                                                       \label{15.28.03}
M_{k,h}(t)\varphi :=&\sum_{|\alpha|\leq1}
b^{\alpha}_k(t)\delta^{\alpha} \varphi ,
\end{align}
for functions $\varphi$ defined on $\mathbb G$. 
It is not difficult to see that 
taking the triple $V_n:=W_{h,2}^1$, $H_n:=W^0_{h,2}$, 
$V_n^{\ast}=(W_{h,2}^1)^{\ast}$, problem 
\eqref{12.28.03}-\eqref{13.28.03} 
can be cast into equation 
\eqref{11.03.04},  and we can easily check that   
Assumption \ref{assumption S}  
and equation \eqref{In}
hold.  Thus \eqref{12.28.03}-\eqref{13.28.03} has a
unique  continuous $W^0_{h,2}$-valued solution 
$v=v^h$ such that for every $h\in (0,1)$, 
$$
E\sup_{t\in[0,T]}|v^h(t)|_{h,0}^2
+E\int_0^T|v^h(t)|_{h,1}^2\,dt \leq M  <\infty. 
$$
Assume now that $d=1$. 
 Consider the normal triple 
$V\hookrightarrow H\equiv H^*\hookrightarrow V^*$ 
with 
 $V:=W^1_2(\mathbb R)$, $H:=W^0_2(\mathbb R)$ 
and $V^*\equiv W^{-1}_2(\mathbb R)$. 
Notice that  
Using \eqref{16.28.03} we can see that 
there is a constant $C$ 
such that almost surely for all $t\in[0,T]$
\begin{eqnarray*}
|D^{\alpha}(a^{\alpha\beta}(t)D^{\beta}\varphi)-
\delta^{\alpha}_{-}(a^{\alpha\beta}(t)
\delta^{\beta}_{+}\varphi)|_{h,0}
&\leq & Ch|\varphi|_{W^3_2(\mathbb R)}, \\
|b^{\alpha}_k(t)D^{\alpha}\varphi-
b^{\alpha}_k(t)\delta^{\alpha}\varphi|_{h,0}
&\leq & Ch|\varphi|_{W^2_2(\mathbb R)}, \\
|F_h(t, D\varphi,\varphi)
-F_h(t, \delta\varphi,\varphi)|_{h,0}
&\leq & C  |h|_{W^2_2(\mathbb R)}  
\end{eqnarray*}
for all $\varphi\in W^3_2(\mathbb R)$ 
and $h\in (0,1)$. Hence the consistency condition 
${\bf (Cn)}$ holds with 
$\mathcal V=W^{3}_2(\mathbb R)
$ 
and $\varepsilon_n=h$. 
Set  $\mathcal H=W^2_2(\mathbb R)$. Assume 
(A1)-(A4) with $l=2$. Then the assumptions of 
Theorem \ref{raten} are satisfied. Thus 
there is a constant $C$ 
such that 
$$
E\sup_{t\in[0,T]}|u(t)-v^h(t)|_{h,0}^2
+E\int_0^T|u(t)-v^h(t)|_{h,1}^2\,dt\leq Ch^2
$$
for all $h\in(0,1)$. 
Now 
we approximate \eqref{13.28.03} by the following Euler 
approximation schemes: 
\begin{eqnarray}
w_{i+1} &=& w_{i}+
\big( L_h(t_{i+1}) w_{i+1} 
+F_h(t_{i+1},\nabla_h w_{i+1},w_{i+1})\big)\,\tau     
\nonumber \\
                                              \label{21.06.06.07}
&&+\sum_{k}\big(M_{k,h}(t_i)w_i 
+g_{k,h}(t_i)\big)(W^k(t_{i+1})-W^k(t_{i})), 
\quad w_0=u_0, \\\
u_{i+1}&=&u_{i}+
\big( L_h(t_{i})u_{i}
+F_h(t_{i+1},\nabla_h u_{i+1}, u_{i} )\big)\tau     
\nonumber                                                      \\
&&+\sum_{k}\big(M_{k,h}(t_i)u_i 
+g_{k,h}(t_i)\big)(W^k(t_{i+1})-W^k(t_{i})), 
\quad v_0=u_0.                                             \nonumber
\end{eqnarray}
for $i=0,1,2,\dots,m-1$, $\tau=T/m$, $t_i=i\tau$. 
Then by Proposition \ref{22.11.04} we get the existence of  
a unique $W^1_{h,2}$-valued solution $w_i$ of 
\eqref{21.06.06.07}, such that  $w_i$ is $\mathcal  F_{t_i}$-
measurable for $i=1,2,\dots,m$, if $\tau$ is sufficiently small. 
By Theorem \ref{th1spi}  
for $e^{h,\tau}_i=(  u(t_i,z)  -w_i(z))_{z\in{\mathbb G}}$,
we get  

$$
E  \max_{0\leq i\leq m}  
|e^{h,\tau}_i|^2_{W^0_{h,2}}+\tau \sum_{1\leq i\leq m}  
E|e^{h,\tau}_i|^2_{W^1_{h,2}}\leq C(\tau+h^2) 
$$
with a constant $C$ independent of $\tau$ and $h$.  
Recall that $\vartheta(n)=\kappa^2/h_n^2$ for any 
sequence $h_n\in(0,1)$ by 
\eqref{23.06.06.07}. Set
$e^{h}_{i,\tau}
=( u(t_i,z) -u_i(z))_{z\in{\mathbb G}}$.  
Then applying Theorem \ref{thexp} we get 
$$
E \max_{0\leq i\leq m}  |e^{h}_{\tau,i}|^2_{W^0_{h,2}}
+\tau  \sum_{0\leq i<m} 
E|e^{h}_{\tau,i}|^2_{W^1_{h,2}}\leq Ch^2,  
$$
with a constant $C$ independent of $\tau$ and $h$, provided 
\eqref{stability} holds with $\kappa^2/h^2$ in place of 
$\vartheta(n)$.
To obtain the corresponding results 
when $d>1$ we need more regularity
in the space variable  from the solution $u$ of
\eqref{ex1eqnl}-\eqref{ini}.  
Assuming  more regularity on the data,  it is
possible  to get the required regularity of $u$. We 
do not want to prove in this paper further 
results on regularity of the solutions to 
 \eqref{ex1eqnl}. Instead of that  we consider the case 
of linear equations, i.e., when $F$ does not depend 
on $p$ and $r$, since in this case the necessary 
results on regularity of the solutions are well 
known in the literature. 
(See  e.g. \cite{KrRo} and \cite{R}.) 

\subsection{Linear stochastic PDEs}

We consider again equation 
\eqref{ex1eqnl}-\eqref{ini} and assume that 
$F=F(t,x,p,r)$ does not depend on $p$ and $r$. 
We fix and integer $l\geq0$. 
Instead of {\bf (A4)}  we assume the following.  

\noindent
{\bf Assumption (A*4)}
{\it
$F(t,x,p,r)=f(t,x)$ and $g_k(t,x)$ are 
  ${\mathcal P}\otimes {\mathcal B}({\mathbb R^d})$
-measurable functions of $(t,\omega,x)$, and their  
derivatives in $x$ up to order $l$  
are 
${\mathcal P}\otimes {\mathcal B}({\mathbb R})$-measurable 
functions such that 
$$
|f(t,.)|_l^2 + \sum_k |g_k(t,.)|_l^2 \leq \eta, 
$$
where $\eta$ is a random variable such that 
 $E\eta \leq M$.  
}

Instead of {\bf(A5)} we make the following assumption. 

\noindent
{\bf Assumption (A*5)} 
{\it Almost surely 

(i)
$
\sum_k|D^{\gamma}(b^{\alpha}_k(t)-b^{\alpha}_k(s))|
\leq K |t-s|^{\frac{1}{2}},
\quad                                                     
\sum_k|g_k(s)-g_k(t)|_l^2 \leq \eta\, |t-s|.
$

(ii)
$
|D^{\gamma}(a^{\alpha,\beta}(t)-a^{\alpha,\beta}(s))|
\leq K|t-s|^{\frac{1}{2}},
\quad
|f(t)-f(s)|^2_l
\leq \eta\, |t-s|
$

\noindent
for all $|\gamma|\leq l$, 
$s,t\in [0,T]$, $x\in {\mathbb R}^d$ and multi-indices 
$|\alpha|\leq1$ and $|\beta|\leq1$,   
where $K$ is a constant and $\eta$ is a random variable 
such that $E\eta\leq M$.}

\smallskip

Consider the space-time discretizations 
with finite differences. The implicit 
and the explicit 
approximations, $v^{h,\tau}$ and 
$v^{h}_{\tau}$ 
are given by the systems of  equations 
defined for $i=0, \cdots, m-1$ by 
\begin{eqnarray}
v^{h,\tau}(t_{i+1})&=&v^{h,\tau}(t_{i})
+\tau \big(L_h(t_{i+1})
v^{h,\tau}(t_{i+1})+f(t_{i+1})\big)       \nonumber \\
                          \label{11.30.03}
&& +\sum_k
\big(M_{k,h}v^{h,\tau}(t_{i})+g(t_{i})\big)
\big(W^k(t_{i+1})-W^k(t_{i})\big), 
                  \\
                                      \label{12.30.03}
v^{h,\tau}(0,z) & =&u(0,z),
\quad z\in{\mathbb G} ,
\end{eqnarray}
and
\begin{eqnarray}
v^{h}_{\tau}(t_{i+1})&=&v^{h}_{\tau}(t_{i})
+\tau (L_h(t_{i})
v^{h}_{\tau}(t_{i})+f(t_{i}))             \nonumber \\
                                          \label{13.30.03}
& & + \sum_k \big(M_{k,h}v^{h}_{\tau}(t_{i})+g(t_{i})\big) 
\big(W^k(t_{i+1})-W^k(t_{i})\big), 
                          \\
                                            \label{14.30.03}
v^{h}_{\tau}(0,z)  &= & u(0,z),
\quad z\in{\mathbb G}, 
\end{eqnarray} 
respectively, where $t_i=i\tau=iT/m$, 
$v^{h,\tau}(t_i)$ and $v^{h}_{\tau}(t_i)$ 
are functions on $\mathbb G$,  
$L_h(t)$ and $M_{k,h}(t)$ are defined 
by \eqref{14.28.03} and \eqref{15.28.03}. 

Take $H_n:=W^0_{h,2}$ and the normal triple 
$
V_n
\hookrightarrow H_n\equiv H_n^*
\hookrightarrow V_n^*
$  
with $V_n:=W^1_{h,2}$. 
Then it is easy to see that 
\begin{equation}                        \label{13.01.04}
(L_h (t_i)\varphi,\psi)_n\leq 
C|\varphi|_{V_n}||\psi|_{V_n}, 
\quad 
(M_{k,h}(t_i)\varphi,\psi)_n\leq 
C|\varphi|_{V_n}||\psi|_{H_n}
\end{equation}
for all $\varphi, \psi\in V_n$, 
where $(\cdot,\cdot)_n$ denotes  
the inner product in $H_n$, 
and $C$ is a constant 
depending only on $d$ and the constant $K$ 
from   Assumption {\bf (A3)}. Thus we can 
define $L_h(t_i)$ and $M_{h,k}(t_i)$ as 
bounded linear operators from $V_n$ into 
$V_n^*$ and $H_n$ respectively. 
 Due to \eqref{13.31.03}    
and \eqref{2.13},  the restriction of 
$u_0$, $f(t_i)$ and $g_k(t_i)$  onto $\mathbb G$ 
are $H_n$-valued random variables 
such that 
$$
E|f(t_i)|_{H_n}^2
\leq p^2 E|f(t_i)|^2_{l}, 
\quad
E|g_k(t_i)|_{H_n}^2
\leq p^2 E|g(t_i)|^2_{l}, 
$$
\begin{equation*}                                   \label{14.01.04}
E|u_0|^2_{H_n}
\leq p^2 E|u_0|^2_{l}, 
\end{equation*}
where $p$ is the constant from \eqref{13.31.03}. 
Moreover, 
\begin{equation}                                   \label{15.01.04}
2(L_h(t_i)\varphi,\varphi)_n
+\sum_{k}|M_{h,k}\varphi|^2_{H_n}
\leq -\tfrac{\lambda}{2}|\varphi|^2_{V_n}
+C|\varphi|^2_{H_n}
\end{equation}
for all $\varphi\in V_n$, where $C$ 
is a constant depending only 
on $d$ and on the constant $K$ from Assumption {\bf (A2)}. 
Thus using the notation 
$u^{n,\tau}_i=v^{h,\tau}(t_i)$, 
$u^{n}_{\tau,i}= v^{h, \tau}(t_i) $
and defining  
\begin{equation*}                                   \label{16.01.04}
A^{n,\tau}_i(\varphi)
=L_h(t_i)\varphi+f(t_i), 
\quad 
B^{n,\tau}_{k,i}(\varphi)
=M_{k,h}(t_i)\varphi+g_k(t_i) 
\end{equation*}
for $\varphi\in W^1_{h,2}$, we can cast 
\eqref{11.30.03}--\eqref{12.30.03} and 
\eqref{13.30.03}--\eqref{14.30.03} into 
\eqref{unspi} and into \eqref{defexp}, 
respectively, and we can see that 
Assumption \ref{assumption ST} and 
condition  \eqref{Inm} hold.  
Consequently, by virtue of
Proposition \ref{22.11.04},   for sufficiently small $\tau$   
\eqref{11.30.03}--\eqref{12.30.03} has a unique 
solution $\{v^{h,\tau}(t_i)\}_{i=0}^m$, 
such that $v^{h,\tau}(t_i)$ is a 
$W^1_{h,2}$-valued 
${\mathcal F}_{t_i}$-measurable 
random variable and  
$E|v^{h,\tau}_i|^2_{h,2}<\infty$.  
Furthermore, by virtue 
of Proposition \ref{10.12.04},  
\eqref{13.30.03}--\eqref{14.30.03} 
has a unique solution 
$\{v^{h}_{\tau}(t_i)\}_{i=0}^m$, such that 
$v^{h}_{\tau}(t_i)$ is a 
$W^1_{h,2}$-valued 
${\mathcal F}_{t_i}$-measurable 
random variable and  
$E|v^{h}_{\tau}(t_i)|^2_{h,2}<\infty$.  

Let $r\geq0$ be an integer, and assume 
that 
 
\begin{equation}                            \label{11.31.03}
l> r+2+\frac{d}{2}.
\end{equation}
Then Theorem \ref{th1spi} gives the following 
result. 
\begin{Th}                                   \label{16.12.04}
Let Assumptions {\bf (A1)}, {\bf (A2)}, {\bf (A3)}, 
{\bf (A*4)} and {\bf (A*5)} hold with $l$ satisfying 
\eqref{11.31.03}. 
Then for 
sufficiently small $\tau$ 
$$
E\max_{1\leq i\leq m}  |v^{h,\tau}(t_i)-u(t_i)|_{h,r}^2+
E \sum_{1\leq i\leq m}  \tau \,
|v^{h,\tau}(t_i)-u(t_i)|_{h,r+1}^2
\leq
C(h^2+\tau ) 
$$
for all $h\in(0,1)$, where $C=C(r,l,p,\lambda, T,K, M,d,d_1)$ 
is a constant. 
\end{Th}
\begin{proof} 
Take  $H_n:=W^r_{h,2}$, 
$H:=W^{l-2}_{2}(\mathbb R^d)$, 
${\mathcal H}:=W^{l}_{2}(\mathbb R^d) $   
and the normal triples 
$$
V_n
\hookrightarrow H_n\equiv H_n^*
\hookrightarrow V_n^*,\quad 
V
\hookrightarrow H\equiv H^*
\hookrightarrow V^*, 
\quad 
{\mathcal V}
\hookrightarrow {\mathcal H}\equiv {\mathcal H^*}
\hookrightarrow {\mathcal V}^{\ast}
$$  
where  
$V_n:=W^{r+1}_{h,2}$, 
$V_n^{\ast}\equiv W^{r-1}_{h,2}$,   
$V:=W^{l-1}_2(\mathbb R^d)$, 
$V^{\ast}\equiv W^{l-3}_2(\mathbb R^d)$, 
${\mathcal V}:=W^{l+1}_{2}(\mathbb R^d)$ and 
${\mathcal V}^{\ast}\equiv W^{l-1}_{2}(\mathbb R^d)=V$. 
Then due to \eqref{11.31.03} there 
is a constant $p$ such that  
for $\Pi_n:=R_h$,  
$$
|\Pi_n\varphi|_{V_n}\leq p|\varphi|_V, 
$$
for all $\varphi\in V$, by virtue of 
\eqref{17.28.03}. It is easy to check 
that 
\eqref{13.01.04}--\eqref{15.01.04} still 
hold, and hence 
\eqref{11.30.03}--\eqref{12.30.03}, written 
as equation \eqref{unspi}, 
satisfies 
Assumption \ref{assumption ST}
and 
condition \eqref{Inm}
 in the new triple 
as well. 
Using \eqref{17.28.03} it is easy  to 
show that due to Assumption {\bf (A3)} 
\begin{eqnarray*}
|L(t_i)\varphi-L_h(t_i)\varphi|_{V^*_n}
&\leq & |L(t_i)\varphi-L_h(t_i)\varphi|_{H_n} 
\leq Ch|\varphi|_{W^{l+1}_2(\mathbb R^d)},   \\
\sum_k
|M_{k}(t_i)\varphi-M_{k,h}(t_i)\varphi|_{H_n}
&\leq & Ch|\varphi|_{W^l_2(\mathbb R^d)}
\end{eqnarray*}
for all $\varphi\in W^{l+1}_2(\mathbb R^d)$, 
where $C$ is a constant depending on $d$, $l$, 
$r$ and on the constant $K$ from Assumption {\bf(A3)}. 
Hence we can see that {\bf (Cn$\tau$)} holds with 
$\varepsilon_n=h$. 
Due to Assumption {\bf (A*5)} we have 
$$
|L(t)\varphi-L(s)\varphi|^2_{V}\leq C|t-s|, 
\quad
|f(t)-f(s)|_{V}^2\leq \eta |t-s|,
$$
$$
\sum_k|M_k(t)\varphi-M_k(s)\varphi|^2_{V}\leq C|t-s|, 
\quad
\sum_k|g_k(t)-g_k(s)|_{V}^2\leq \eta |t-s|, 
$$
where $\eta$ is the random variable from Assumption {\bf (A*5)}, 
and $C$ is a constant depending on 
$d$, $d_1$, $l$ and on the constant 
$K$ from Assumption {\bf (A*5)}. It is an easy 
exercise to show that due to Assumptions 
{\bf (A3)} and {\bf (A*4)} condition {\bf (R3)} 
from Assumption \ref{assumption R} holds.  
From \cite{KrRo} it is known that under 
the Assumptions 
{\bf (A1)}--{\bf (A3)} and {\bf (A*4)} 
the problem 
\eqref{12.28.03}--\eqref{13.28.03} 
has a unique solution  $u$  
on $[0,T]$, and that $u$ is a continuous 
$W^l_2(\mathbb R^d)$-valued $(\mathcal F_t)$-adapted 
stochastic process such that 
\begin{align*}
E\sup_{t\in[0,T]} & |u(t)|^2_{l}  +
E\int_0^T  |u(t)|^2_{l+1}  \,dt\leq 
\\
&  CE|u_{0}|^2_{l} 
+CE\int_0^T\Big(|f(t)|^2_{{l-1}}
+\sum_k|g_k(t)|^2_{l}\Big)\,dt, 
\end{align*}
where $C$ is a constant depending on $d$, $d_1$ and 
the constants $\lambda$ and $K$ from 
Assumptions {\bf (A1)}, {\bf (A3)} and {\bf (A*4)}. 
Hence the regularity   conditions  
 {\bf (R1)} and {\bf (R2)} 
in Assumption \ref{assumption R}
clearly hold. Now we  can conclude the proof by 
applying Theorem \ref{th1spi}. 
\end{proof}

Let us now investigate the rate of convergence 
of the explicit space-time approximations. 
Take the normal triple 
$V_n
\hookrightarrow H_n\equiv H_n^*
\hookrightarrow V_n^*
$ 
with $V_n:=W^{r+1}_{h,2}$, $H_n:=W^{r}_{h,2}$, 
and notice that due to Assumption {\bf (A3)}
\begin{equation}                        \label{13.11.04}
(L(t_i)\varphi,\psi)_n\leq 
C_1|\varphi|_{V_n}|\psi|_{V_n}, 
\quad 
(M_{k,h}(t_i)\varphi,\psi)_n\leq 
C_{2k}|\varphi|_{V_n}|\psi|_{H_n}
\end{equation}
with some constants $C_1$ and $C_{2k}$ depending 
only on $d$, $r$ and the constant $K$ 
from Assumption {\bf (A3)}.  Set $L_1=C_1^2$ 
and $L_2=\sum_kC_{2k}^2$. 
Then Theorem \ref{thexp} yields 
the following theorem, which improves a result 
from \cite{HY}. 
\begin{Th}                               \label{14.11.04}
Let Assumptions {\bf (A1)}, {\bf (A2)}, {\bf (A3)}, 
{\bf (A*4)} and {\bf (A*5)} hold with $l$ satisfying 
\eqref{11.31.03}. 
Let $h$ and $\tau$ satisfy 
\begin{equation}                         \label{stability2} 
L_1\kappa^2\frac{\tau}{h^2} 
+2\kappa(L_1L_2)^{1/2}\frac{\sqrt{\tau}}{h}\leq q
\end{equation}
for a constant $q<\lambda$.  Then 
$$
E\max_{1\leq i\leq  m}|v^h_{\tau}(t_i)-u(t_i)|_{h,r}^2+
E\sum_{0\leq i< m}  \tau \, |v^h_{\tau}(t_i)-u(t_i)|_{h,r+1}^2
\leq
C(h^2+\tau ) 
$$
for all $h\in(0,1)$, where $C=C(r,l,p,\lambda, q, T,K, M,d,d_1)$ 
is a constant. 
\end{Th}
\begin{proof} As in the proof of 
Theorem \ref{16.12.04} we  
take $H_n:=W^r_{h,2}$, 
$H:=W^{l-2}_{2}(\mathbb R^d)$, 
${\mathcal H}:=W^{l}_{2}(\mathbb R^d)$   
and the normal triples 
$$
V_n
\hookrightarrow H_n\equiv H_n^*
\hookrightarrow V_n^*,\quad 
V
\hookrightarrow H\equiv H^*
\hookrightarrow V^*, 
\quad 
{\mathcal V}
\hookrightarrow {\mathcal H}\equiv {\mathcal H}^*
\hookrightarrow {\mathcal V}^{\ast}
$$
with $V_n:=W^{r+1}_{h,2} 
$, $V:=W^{l-1}_{2}(\mathbb R^d)$, 
${\mathcal V}=W^{l+1}_{2}(\mathbb R^d)$,  we cast 
\eqref{13.30.03}--\eqref{14.30.03} into \eqref{defexp}, 
and see that Assumptions \ref{assumption R} and 
\ref{assumption ST}, conditions  {\bf (Cn$\tau$)} and 
\eqref{Inm}  
 of Theorem \ref{thexp} hold.  
Furthermore, $\vartheta(n)=\frac{\kappa^2}{h^2}$. 
We can easily check 
that by virtue of \eqref{13.11.04} and \eqref{15.12.04}, 
condition \eqref{stability2} yields condition 
\eqref{stability}. Hence applying Theorem \ref{thexp} 
we finish the proof. 
\end{proof}

\begin{corollary} Let  $k\geq 0$ be an integer  and   
let Assumptions {\bf (A1)}, {\bf (A2)}, {\bf (A3)}, 
{\bf (A*4)} and {\bf (A*5)} hold with $l$ satisfying 
  $ l> k+2+d $.   
Then the following statements are valid 
for all multi-indices $|\alpha|\leq k$:  

(i) For 
sufficiently small $\tau$ 
$$
E\max_{1\leq i\leq m}\sup_{z\in\mathbb G}
|\delta^{\alpha}(v^{h,\tau}(t_i,z)-u(t_i,z))|
\leq
C(h+\sqrt\tau )
$$
holds for all $h\in(0,1)$, where
$ C=C(l,p,\lambda, T,K,M,d,d_1)$ 
is a constant.

(ii) Assume also that $\tau$ and $h$ 
satisfy \eqref{stability2}. 
Then  $$
E\max_{1\leq i\leq m}\sup_{z\in{\mathbb G}}
|\delta^{\alpha}(v^h_{\tau}(t_i,z)-u(t_i,z))|
\leq
C(h+\sqrt{\tau })
\leq C\left(1+\kappa^{-1}\sqrt{\lambda/L_1}\right)h
$$
for all $h\in(0,1)$, where $C=C(r,l,p,\lambda, q, T,K, M,d,d_1)$ 
is a constant. 
\begin{proof} By the discrete version of Sobolev's theorem on embedding 
$W^2_m(\mathbb R^d)$ into ${\mathcal C}^k(\mathbb R^d)$ one knows 
that if $m\geq k+\tfrac{d}{2}$, then 
$$
\sup_{z\in\mathbb G}|\delta^{\alpha}\varphi(z)|
\leq C|\varphi|_{W^m_{h,2}}
$$
for all $h\in(0,1)$, 
$\varphi\in W^m_{h,2}$ and $|\alpha|\leq k$, 
 where $C=C(d,m,k)$ is a constant  
(see e.g. \cite{HY}). Hence the above statements follow 
immediately from 
Theorems \ref{16.12.04} and \ref{14.11.04}. 
\end{proof}
\end{corollary}

\noindent{\bf Acknowledgments:} 
The authors thank 
A. Cohen, G. Kerkyacharian and G. Gr\"un 
for stimulating discussions
about wavelets and  finite elements.  
They are also grateful to the referee 
for his comments and suggestions 
which helped to improve the
presentation of the paper.

\noindent {\it Addresses :}\\
Istvan Gy\"ongy,\\
 School of Mathematics
University of Edinburgh,\\
{\it and} Maxwell Institute for Mathematical Sciences,\\
 King's  Buildings,\\
 Edinburgh, EH9 3JZ, United Kingdom\\
email : gyongy@maths.ed.ac.uk
\smallskip

and
\smallskip

\noindent  Annie Millet \\
Laboratoire de Probabilit\'es et Mod\`eles Al\'eatoires\\
  Universit\'es Paris~6-Paris~7, Bo\^{\i}te Courrier 188,\\
      4 place Jussieu, 75252 Paris Cedex 05, France\\
      {\it  and } SAMOS-MATISSE,
      Centre d'\'Economie de la Sorbonne,\\
 Universit\'e Paris 1 Panth\'eon Sorbonne,\\
       90 Rue de Tolbiac,\\ 75634 Paris Cedex 13, France\\
email : annie.millet@upmc.fr {\it and} amillet@univ-paris1.fr
\smallskip

\end{document}